\documentclass[12pt]{amsart}
\usepackage{txfonts}
\usepackage{amsfonts}
\usepackage{amsfonts}
\usepackage{amssymb,latexsym}
\usepackage{cite}
\usepackage{enumerate}
\newtheorem{theorem}{Theorem}

\newtheorem{proposition}[theorem]{Proposition}
\newtheorem{lemma}[theorem]{Lemma}
\newtheorem{definition}[theorem]{Definition}

\newcommand {\p}{\partial}

\numberwithin{equation}{section}
\numberwithin{theorem}{section}

\usepackage{titletoc}

\usepackage{hyperref}
\hypersetup{hypertex=true,colorlinks=true,linkcolor=blue,anchorcolor=blue,citecolor=blue}
\begin{document}
\title[The homogeneous complex k-Hessian equation]{The exterior Dirichlet problem for homogeneous complex $k$-Hessian equation}
\author{Zhenghuan Gao}
\address{Department of Mathematics, Shanghai University, Shanghai, 200444, China}
\email{gzh@shu.edu.cn}
\author{Xinan Ma} 
\address{School of Mathematical Sciences, University of Science and Technology of China, Hefei 230026,
	Anhui Province, China}
\email{xinan@ustc.edu.cn}

\author{Dekai Zhang}  
\address{Department of Mathematics, Shanghai University, Shanghai, 200444, China}
\email{dkzhang@shu.edu.cn}

\begin{abstract}
 In this paper, we consider the homogeneous  complex $k$-Hessian equation in an exterior domain $\mathbb {C}^n\setminus \Omega$. 
 We prove the existence and uniqueness  of the  $C^{1,1}$ solution by constructing approximating solutions.
 The key point for us is to establish the uniform gradient estimate and the second order estimate. 

\end{abstract}
\maketitle
\tableofcontents 

\section{Introduction}
Let $u$ be a real $C^2$ function in $\mathbb C^n$ and $\lambda=(\lambda_1,\cdots,\lambda_n)$ be the eigenvalues of the complex Hessian $(\frac{\partial^2u}{\partial z_j\partial \bar z_k})$, the complex $k$-Hessian operator is defined by
\begin{align}
	H_{k}(u):=\sum\limits_{1\le i_1<\cdots i_k\le n}\lambda_{i_1}\cdots \lambda_{i_k},
\end{align}
where $1\le k\le n$. Using the operators $d=\partial+\overline\partial$ and $d^c=\sqrt{-1}(\overline \partial-\partial )$, such that $dd^c=2\sqrt{-1}\partial\overline \partial$, one gets
$$(dd^cu)^k\wedge\omega^{n-k}=4^nk!(n-k)!H_k(u)d\lambda,$$
where $\omega=dd^c|z|^2$ is the fundamental K\"{a}hler form and $d\lambda$ is the volume form.
When $k=1$, $H_1(u)=\frac14\Delta u$. When $k=n$, $H_n(u)=\det u_{i\bar j}$ is the complex Monge-Amp\`ere operator. 

Let $\Omega$ be a bounded smooth domain in $\mathbb{C}^n$, the Dirichlet problem for the complex $k$-Hessian equation is as follows
\begin{align}\label{khessian}
	\left\{\begin{aligned}
		H_k(u) =&f \qquad\text{in} \quad \Omega,\\
		u=&\varphi \qquad\text{on} \quad \partial\Omega,
		\end{aligned}
	\right.
\end{align}
where $f$ and $\varphi$ are given smooth functions.
When $k=1$, the $k$-Hessian equation is the Poisson equation.
When $k=n$, it is the well known  complex Monge-Amp\`ere equation.
\subsection{Some previous results}
We briefly give some studies on the Dirichlet problem for the $k$-Hessian equation and the complex $k$-Hessian equation in the  nondegenerate case i.e. $f>0$ and in the degenerate cases i.e. $f\ge 0$. In general, the $k$-Hessian equation (the complex $k$-Hessian equation) is a fully nonlinear equation.
\subsubsection{Results on bounded domains }
For the $k$-Hessian equation in $\mathbb R^n$, if $f>0$,
Caffarelli-Nirenberg-Spruck \cite{CNSIII} solved the Dirichlet problem in a bounded $(k-1)$-convex domain.  Guan \cite{guan1994cpde} solved the Dirichlet problem by only assuming the existence of a subsolution.  For the complex $k$-Hessian equation in $\mathbb C^n$, Li \cite{Lisongying2004} solved \eqref{khessian} in a bounded $(k-1)$-pseudoconvex domain.

 There are lots of studies on the Dirichlet problem in bounded domains in $\mathbb R^n$ of  degenerate fully nonlinear equations.
Caffarelli-Nirenberg-Spruck \cite{cns1986rmi} show the $C^{1,1}$ regularity of the homogeneous Monge-Amp\`ere equation  i.e. $f\equiv 0$.
If  $f^{\frac{1}{n-1}}\in C^{1,1}$, Guan-Trudinger-Wang \cite{guantrudingerwang1999acta} proved the optimal $C^{1,1}$ regularity result due to the counterexample by Wang \cite{wang1995pams}.
The  $k$-Hessian equation case was proved by Krylov\cite{krylov1989, krylov1994} and Ivochina-Trudinger-Wang\cite{itw2004cpde} (PDE's proof) by assuming $f^{\frac{1}{k}}\in C^{1,1}$.
Dong \cite{dong2006cpde} proved the $C^{1,1}$ regularity  for some  degenerate mixed type Hessian equations.

For the Dirichlet problem of degenerate complex Monge-Amp\`ere equation, Lempert \cite{Lempert1983} showed that $(dd^cu)^n=0$ in a punched strictly convex domain $\Omega\backslash\{z\}$ with logarithm growth near $z$  admits a unique real analytic solution.
Zeriahi \cite{Zeriahi2013} studied the viscosity solution to the Dirichlet problem of degenerate complex Monge-Amp\`ere equation.  

\subsubsection{Results on unbounded domains }
There are lots of results on the exterior Dirichlet problem for viscosity solutions of nondegenerate  fully nonlinear  equations. The $C^0$ viscosity solution  for the  Monge-Amp\`ere equation: $\det D^2 u=1$ with  prescribed asymptotic behavior at infinity was obtained by Caffarelli-Li\cite{cl2003cpam}.
 The  $k$-Hessian equation case  was showed by Bao-Li-Li \cite{baolili2014tams}. For the related results  on  other type nondegenerate fully nonlinear equations, one can see \cite{baoli2013, lili2018jde, libao2014jde, li2019tams}. 

 The global $C^{k+2,\alpha}$ regularity of the homogeneous Monge-Amp\`ere equation
 in a strip region was proved by Li-Wang\cite{liwang2015dcds}
 by assuming that the boundary functions are locally uniformly convex and $C^{k,\alpha}$. They showed that the uniform convexity of the boundary functions is necessary.
 
 For $1\le k<\frac{n}{2}$, the $C^{1,1}$ regularity of Dirichlet problem for the homogeneous $k$-Hessian equation in $\mathbb R^n\setminus \overline \Omega$ was proved by Xiao \cite{xiao2022} by assuming that the domain $\Omega$ is $(k-1)$-convex and starshaped.
 For $1\le k\le n$, Ma-Zhang \cite{MaZhang2022} proved the $C^{1,1}$ regularity  of the $k$-Hessian equation when $\Omega$ is convex and strictly $(k-1)$ convex.  The prescribed asymptotic behavior is  $\log |x|+O(1)$ if $k=\frac{n}{2}$ and $|x|^{2-\frac{n}{k}}+O(1)$ if $k>\frac{n}{2}$.
\subsection{Motivation}
Our research is motivated by the study of regularity of extremal  function. 
For the smoothly strictly convex domain $\Omega$, Lempert \cite{lempert1985duke} prove the pluricomplex Green function in $\mathbb C^n\setminus \Omega$ is smooth (analytic). 
In \cite{gpf2002am,gpf2010}, P. F. Guan proved the $C^{1,1}$ regularity of the solution to the homogeneous complex Monge-Amp\`ere equation in a ring domain.  
Then he solved a conjecture of Chern-Levine-Nirenberg on the extended intrinsic norms. For the smoothly strongly pseudoconvex domain $\Omega$, 
B. Guan \cite{gb2007imrn} proved the $C^{1,1}$ regularity and decay estimates of pluricomplex Green function in $\mathbb C^n\backslash \Omega$ by considering the exterior Dirichlet problem for the homogeneous complex Monge-Amp\`ere
equation. 

Another motivation is on the proof of geometric inequalities by considering the
%
%
%
When $\Omega$ is $(k-1)$-convex and starshaped,
Guan-Li\cite{guanli2009adv} proved Alexandrov-Fenchel inequalities 
by the inverse curvature flows. If $\Omega$ is $k$-convex, Chang-Wang\cite{cw2013adv}, Qiu\cite{q2015ccm} proved Alexandrov-Fenchel inequalities  by the optimal transport method. Whether Alexandrov-Fenchel inequalities hold for  $(k-1)$-convex domain is still open.
Recently, Agostiniani-Mazzieri \cite{AM2020CVPDE} proved some geometric inequalities such as Willmore inequality by considering the exterior Dirichlet problem of the Laplace equation. 
Fogagnolo and Mazzieri and Pinamonti \cite{fmp2019ihp} showed the volumetric Minkowski inequality by considering the the exterior Dirichlet problem of the $p$-Laplacian equation.
Agostiniani-Fogagnolo-Mazzieri  \cite{afm2022arma} removed the convexity assumption in \cite{fmp2019ihp} for the domain. 
\subsection{Our main result}
In this paper, we consider the following exterior Dirichlet problem for the complex $k$-Hessian equation. 

For $1\le k<n$, since the Green function in this case is $-|z|^{2-\frac{2n}k}$, we consider the $k$-Hessian equation  as follows
\begin{equation}\label{case1Equa1.1}
\begin{cases}
(dd^cu)^k\wedge \omega^{n-k}=0&\quad\text{in }\Omega^c:=\mathbb C^n\backslash \overline\Omega,\\
u=-1&\quad\text{on }\partial\Omega,\\
u(z)\rightarrow 0&\quad\text{as }|z|\rightarrow\infty.
\end{cases}
\end{equation}

\begin{theorem}\label{main07201}
Assume $1\le k<n$. Let $\Omega$ be a smoothly strongly pseudoconvex domain in $\mathbb{C}^n$ such that $0\in \Omega$ and $\overline\Omega$ is holomorphically convex in ball centered at $0$. There exists a unique  $k$-subharmonic solution $u\in C^{1,1}(\overline{\Omega^c})$ of the equation \eqref{case1Equa1.1}. Moreover, there exists uniform constant $C$ such that for any $z\in \overline{\Omega^c}$ the following holds
\begin{align}\label{decay10720}
\begin{cases}
&C^{-1}|z|^{2-\frac{2n}k}\le -u(z)\le C|z|^{2-\frac{2n}{k}},\\
&|Du|(z)\le C|z|^{1-\frac{2n}{k}},\\
&\Delta u(z)\le C|z|^{-\frac{2n}{k}},\\
&|Du|_{C^{0,1}(\Omega^c)}\le C.
\end{cases}
\end{align}

\end{theorem}

Here the $k$-subharmonic function is defined in Section 2 and we use the notation $\overline {\Omega^c}:=\mathbb C^n\setminus \Omega$.
 Let $r_0$ be the constant such that $B_{r_0}\subset\subset\Omega$ and   $R_0, S_0$ be constants such that $\overline\Omega$ is holomorphically convex in $B_{S_0}$ and $\Omega\subset\subset B_{R_0}\subset \subset B_{S_0}$, where $B_{r_0}$, $B_{R_0}$ and $B_{S_0}$ are balls centered at $0$ with radius $r_0$, $R_0$ and $S_0$ respectively.

To prove Theorem \ref{main07201}, we  consider the following approximating equation
\begin{equation*}
\begin{cases}
H_k(u^\varepsilon)=f^\varepsilon\quad\ \text{in }\Omega^c,\\
u^\varepsilon=-1\quad\quad\ \ \  \text{on }\partial\Omega,\\
u^\varepsilon(z)\rightarrow 0\quad\quad\text{as }|z|\rightarrow \infty.
\end{cases}
\end{equation*}
where $f^{\varepsilon}=c_{n,k}\varepsilon^2(1+\varepsilon^2)^{n-k}(|x|^2+\varepsilon^2)^{-n-1}$ (see Section \ref{section4}). 

$u^{\varepsilon}$ will be obtained by approximating solutions $u^{\varepsilon,R}$ defined on bounded domains: $\Sigma_R:=B_R\setminus\overline\Omega$ (see Section \ref{section4}).
The existence and uniqueness of the smooth $k$-subharmonic solution of $u^{\varepsilon,R}$ follows from Li \cite{Lisongying2004} if we can construct a subsolution. 
The key point is to establish the uniform $C^{2}$ estimates for $u^{\varepsilon,R}$.


In Section \ref{section2}, we give some preliminaries. In Section \ref{section3}, we solve the Dirichlet problem of degenerate complex $k$-Hessian equation in a ring domain. 
Section \ref{section4} is the main part of this paper. We show uniform $C^{1,1}$ estimate of the solution which is the limit of the solutions of nondegenerate complex $k$-Hessian equation. The key ingredient is to establish  uniform gradient estimates and uniform second order estimates. We use the idea of Hou-Ma-Wu \cite{HMW} (see also Chou-Wang\cite{cw2001cpam}) to establish the uniform second order estimates.
Theorem \ref{main07201}
will be proved in Section \ref{section5}. 

\begin{align*}
	\end{align*}

%

\section{Preliminaries}\label{section2}
\subsection{$k$-subharmonic solutions}
In this section we give the definition of $k$-subharmonic functions and definition of $k$-subharmonic solutions.

The $\Gamma_k$-cone is defined by
\begin{align}
	\Gamma_{k}:=\{\lambda\in \mathbb{R}^n|S_{i}(\lambda)>0, 1\le i\le k\}
\end{align}
Recall $S_k(\lambda):=\sum\limits_{1\le i_1< \cdots <i_k\le n}\lambda_{i_1}\cdots\lambda_{i_k}$, and $S_k(A):=\delta_{i_1\cdots i_k}^{j_1\cdots j_k}A_{i_1j_1}\cdots A_{i_kj_k}$, where $\delta_{i_1\cdots i_k}^{j_1\cdots j_k}$ is the Kronecker symbol, which has the value $+1$ (respectively, $-1$) if $i_1,i_2\cdots i_k$
are distinct and $(j_1j_2\cdots j_k)$ is an even permutation (respectively, an odd permutation) of
$(i_1i_2\cdots i_k)$, and has the value 0 in any other cases. We use the convention that $S_0(A)=1$. It is clear that 
$S_k(A)=S_k(\lambda(A))$, where $\lambda(A)$ are the eigenvalues of $A$.

One can find the concavity property of $S_k^{\frac{1}{k}}$  in \cite{CNSIII}.
\begin{lemma}\label{concavity}
	$S_k^{\frac{1}{k}}$ is a concave function in $\Gamma_k$. In particular, $\log S_k$ is concave in $\Gamma_k$.
	\end{lemma}
The following facts about elementary symmetric polynomial  are useful in proving gradient estimates.	
\begin{proposition}\label{2ineqs}
We have the following two inequalities,
\begin{itemize}
\item[(a)] If $\lambda \in\Gamma_k,$ then
$$\frac{S_k^2(\lambda| i)}{S_{k-1}(\lambda| i)}\geq \frac{k+1}k\frac{n-k}{n-k-1}S_{k+1}(\lambda| i);$$
\item[(b)] If $\lambda \in\Gamma_k,$ then
$$\frac{S_{k}(\lambda| i)}{S_{k-1}(\lambda| i)}\leq \frac1k\frac{n-k}{n-1}S_1(\lambda| i).$$
\end{itemize}
\end{proposition}
\begin{proof}
Since $\lambda\in\Gamma_k$, we have $S_{k-1}(\lambda| i)>0$. The first inequality follows from Newton inequality. Now we prove (b). Since $\lambda\in \Gamma_k$, we have $S_h(\lambda| i)>0$, $\forall\, h=0,1,\cdots,k-1$. If $S_k(\lambda|i)\leq 0$, (b) holds naturally. If $S_k(\lambda| i)> 0$, the second inequality follows from the generalized Newton-MacLaurin inequality.
\end{proof}
The following two propositions enable us to adopt a casewise argument to deal with the third order terms as in \cite{cw2001cpam} and \cite{HMW}.
\begin{proposition}
Let $\lambda=(\lambda_1,\cdots,\lambda_n)\in\overline \Gamma_k$, and $\lambda_1\geq\lambda_2\geq \cdots\geq\lambda_n$.
Then there exists $\theta=\theta(n,k)>0$ such that 
$$S_{k-1}(\lambda|k)\geq \theta\lambda_1 S_{k-2}(\lambda|1k),$$ from which it follows
\begin{align}\label{eq719::1}S_{k-1}(\lambda| i)\geq \theta \lambda_1\lambda_2\cdots\lambda_{k-1},\quad\forall\, i\geq k.\end{align}
\end{proposition}
The following proposition was proven in \cite{cw2001cpam}. In \cite{HMW}, Hou-Ma-Wu provided a sharp constant $\theta=\frac kn$ in \eqref{eq719::2}.
\begin{proposition}\label{lemmaChouWang}
Let $\lambda=(\lambda_1,\cdots,\lambda_n)\in\overline \Gamma_k$, and $\lambda_1\geq\lambda_2\geq \cdots\geq\lambda_n$.
Then there exists $\theta=\theta(n,k)>0$ such that 
\begin{align}\label{eq719::2}\lambda_1 S_{k-1}(\lambda|i)\geq \theta S_k(\lambda).\end{align}
Moreover, for any $\delta \in(0,1)$ there exists $K>0$ such that if
$$S_k(\lambda)\leq K\lambda_1^k\quad\text{or}\quad |\lambda_i|\leq K\lambda_1\quad\text{for any }i=k+1,k+2,\cdots,n,$$
we have 
\begin{align}\label{eq719::3}\lambda_1 S_{k-1}(\lambda|1)\geq(1-\delta)S_k(\lambda).\end{align}
\end{proposition}

One can see the Lecture notes by Wang \cite{wang2009} for more properties of the $k$-Hessian operator, and see Blocki \cite{BlockiAIF} for those of the complex $k$-Hessian operator.
We following the definition by Blocki \cite{BlockiAIF} to give the definition of $k$-subharmonic functions.
\begin{definition}
Let $\alpha$ be a real $(1,1)$-form in $U$, a domain of $\mathbb C^n$. We say that $\alpha$ is $k$-positive in $U$ if the following inequalities hold
$$\alpha^j\wedge \omega^{n-j}\geq 0, \forall\,j=1,\cdots,k.$$
\end{definition}

\begin{definition}
	Let $U$ be a domain in $\mathbb C^n$. 
	
(1). A function $u:U\rightarrow \mathbb R\cup \{-\infty\}$ is called $k$-subharmonic if it is subharmonic and for all $k$-positive real $(1,1)$-form $\alpha_1,\cdots,\alpha_{k-1}$ in $U$,
$$dd^cu\wedge \alpha_1\wedge \cdots \wedge \alpha_{k-1}\wedge \omega^{n-k}\geq 0.$$
The class of all $k$-subharmonic functions in $U$ will be denoted by $\mathcal{SH}_k(U)$.

(2). A function $u\in C^2 (U)$ is called  \emph{$k$-subharmonic} (\emph{strictly $k$-subharmonic}) if $\lambda(\partial\overline\partial u)\in \overline \Gamma_{k}$ \ ($\lambda$ ($\partial\overline\partial u)\in\Gamma_k$).

%
\end{definition}
If $u\in \mathcal{SH}_{k}(U)\cap C(U)$, $(dd^{c}u)^{k}\wedge \omega^{n-k}$ is well defined in pluripotential theory by Blocki \cite{BlockiAIF}. We need the following comparison principle by Blocki\cite{BlockiAIF} to prove the uniqueness of the continuous solution of the problem \eqref{case1Equa1.1}.
\begin{lemma}\label{comparison0718}
Let $U$ be a bounded domain in $\mathbb C^n$,  $u, v\in \mathcal{SH}_k(U)\cap C(\overline U)$ satisfy
	\begin{align}
		\left\{\begin{aligned}
			(dd^cu)^k\wedge \omega^{n-k}\ge& (dd^cv)^k\wedge \omega^{n-k}\quad \text{in}\ U,\\
			u\le& v\ \quad\qquad\qquad\ \ \ \ \text{on}\ \p U.
			\end{aligned}
		\right.
		\end{align}
Then $u\le v$ in $U$.
\end{lemma}

\subsection{The existence of the subsolution}
\begin{definition}
$\rho$ is called a defining function of $C^1$ domain $U$, if $U=\{z:\rho(z)<0\}$ and $|D\rho|\neq 0$ on $\partial U$.
\end{definition}
\begin{definition}
A $C^2$ domain $U$ is called pseudoconvex (strictly pseudoconvex)  if it is Levi pseudoconvex (strictly Levi pseudoconvex). That is, for a $C^2$ defining function of $U$ defined in a neighborhood of $ U$, 
the Levi form at every point $z\in \p U$ defined by
$$L_{\partial U,z}(\xi)=\frac{1}{|D\rho(z)|}\sum_{j,k}\frac{\partial^2\rho}{\partial z_j\partial\bar z_k}\xi_j\bar\xi_k, \quad\xi\in \ ^hT_{\partial U,z}$$	 is nonnegative (positive). $ ^hT_{\partial U,z}:=\{\xi\in\mathbb C^n\mid \sum_j\frac{\partial\rho}{\partial z_j}\xi_j=0\}$ is the holomorphic tangent space 	to $\partial U$ at $z$.
\end{definition}
\begin{definition}
A $C^2$ domain $U$ is called $k$-pseudoconvex (strictly $k$-pseudoconvex)  if for a $C^2$ defining function of $U$ defined in a neighborhood of $ U$, 
$$\lambda\bigg\{\frac{\partial^2\rho}{\partial z_i\partial\bar z_j}\bigg\}_{1\leq i,j\leq n-1}\in\overline\Gamma_k\ (\in \Gamma_k ),\quad \forall\,z\in\partial U,$$	 where $(z_1,\cdots,z_{n-1})$ is a holomorphic coordinate system of $  ^hT_{\partial U,z}$ near $z $.
\end{definition}

We need the following  lemmas by Guan\cite{gpf2002am} to construct the subsolution of the $k$-Hessian equation in a ring.
\begin{lemma}\label{Guan2002}
Suppose that $U$ is a bounded smooth domain in $\mathbb{C}^n$. For $h, g\in C^m(U)$, $m\ge 2$, for all $\delta>0$, there is an $H\in C^m(U)$ such that

\begin{enumerate}[(1)]
\item
 $H\ge \max\{h, g\}$ \ \text{and}

	\begin{align*}
 H(z)=\left\{ {\begin{array}{*{20}c}
		{h(z), \quad   \text{if } \ h(z)-g(z)>\delta  }, \\
		g(z) , \ \quad \text{if } \ g(z)-h(z)>\delta;\\
\end{array}} \right.
\end{align*}
\item
{There exists}  $|t(z)|\le 1 $ {such that}

\begin{align*}
\left\{H_{i\bar j}(z)\right\}\ge
 \left\{\frac{1+t(z)}{2}g_{i\bar j}+\frac{1-t(z)}{2}h_{i\bar j}\right\},\ \text{for all} \ z\in\left\{|g-h|<\delta\right\}.
 \end{align*}
\end{enumerate}
	\end{lemma}
By Lemma \ref{concavity}, we see $H$ is  $k$-subharmonic if $h$ and $g$ are both $k$-subharmonic.

The following lemma was proved by Guan \cite{gpf2002am}. 
\begin{lemma}\label{subu0720}
	Let $\Omega_0$ and $\Omega_1$ be smooth, strongly pseudoconvex domain in $\mathbb R^n$ with $\Omega_1\subset\subset\Omega_0$ . Assume that $\Omega_1$ is holomorphically convex in $\Omega_0$. Then there exists a strictly plurisubharmonic function $\underline{u}\in C^{\infty}(\overline \Omega)$ with $\Omega:=\Omega_0\setminus\overline\Omega_1$ satisfying
\begin{equation}\label{underlineu071911}
\begin{cases}
H_k(\underline u)\geq \epsilon_0,&\quad\text{in }\Omega,\\
\underline u=\mathrm\tau\rho_1,&\quad\text{near }\partial\Omega_1,\\
\underline u=1+\mathrm K\rho_0,&\quad\text{near }\partial\Omega_0,\\
\end{cases}
\end{equation}
where $\rho_{0}$ and $\rho_1$ are  defining functions of $\Omega_0$ and $\Omega_1$, $\tau$ and $K$ are uniform constants.
\end{lemma}

In \cite{gpf2002am}, Guan considered  the Dirichlet problem of homogeneous complex Monge-Amp\`ere equation in a smooth ring.
\begin{equation}\label{nEquaRing}
\begin{cases}
(dd^cu)^n=0&\quad\text{in }\Omega:=\Omega_0\backslash\overline\Omega_1,\\
u=0&\quad\text{on }\partial\Omega_1,\\
u=1&\quad\text{on }\partial\Omega_0.
\end{cases}
\end{equation}
Guan \cite{gpf2002am} proved the following.
\begin{theorem}\label{ring0720}
	Let $\Omega_0, \Omega_1$ be smooth, strongly pseudoconvex domains and assume that $\Omega_1$ is holomorphically convex in $\Omega_0$.
There exists a unique solution $u\in C^{1,1}(\overline \Omega)$ of the equation \eqref{nEquaRing}.
\end{theorem}

\section{The Dirichlet problem for the homogeneous $k$-Hessian equations in the ring}\label{section3}
In this section, we consider the Dirichlet problem of the homogeneous complex $k$-Hessian equation in a smooth ring.
\begin{equation}\label{EquaRing}
\begin{cases}
(dd^cu)^k\wedge\omega^{n-k}=0,&\quad\text{in }\Omega:=\Omega_0\backslash\overline\Omega_1,\\
u=0,&\quad\text{on }\partial\Omega_1,\\
u=1,&\quad\text{on }\partial\Omega_0.
\end{cases}
\end{equation}
We assume that $\Omega_1\subset\subset\Omega_0$ are smooth, strongly pseudoconvex domains and $\Omega_1$ is holomorphically convex in $\Omega_0$.
Using Lemma \ref{subu0720}, there exists a smooth, strictly plurisubharmonic subsolution $\underline u$ satisfying
\begin{equation}\label{underlineu071911}
\begin{cases}
H_k(\underline u)\geq \epsilon_0,&\quad\text{in }\Omega,\\
\underline u=\mathrm\tau\rho_1,&\quad\text{near }\partial\Omega_1,\\
\underline u=1+\mathrm K\rho_0,&\quad\text{near }\partial\Omega_0,\\
\end{cases}
\end{equation}
where $\mathrm{{\tau}},\mathrm{{K}}$ are positive constants and $\rho_i$ are  defining functions of $\Omega_i$.
\begin{theorem}\label{ring0720}
	Let $\Omega_0, \Omega_1$ be smooth, strongly pseudoconvex domains and assume that $\Omega_1$ is holomorphically convex in $\Omega_0$.
There exists a unique solution $u\in C^{1,1}(\overline \Omega)$ of the equation \eqref{EquaRing}.
\end{theorem}

The uniqueness follows from Lemma \ref{comparison0718}, the comparison principle for $k$-subharmonic solutions to complex $k$-Hessian equations. Next, we prove the existence and regularity of $k$-subharmonic solution by approximation. Indeed, for every $0<\epsilon<\epsilon_0$, we consider the following problem
\begin{align}\label{ApprEquaRing3}
\begin{cases}
H_k(u^\epsilon)=\epsilon&\quad\text{in }\Omega,\\
u^\epsilon=0&\quad\text{on }\partial\Omega_1,\\
u^\epsilon=1&\quad\text{on }\partial\Omega_0.
\end{cases}
\end{align}

Since $\underline u$ in \eqref{underlineu071911} is a subsolution to \eqref{ApprEquaRing3}, by Li \cite{Lisongying2004}, the above problem has a unique smooth solution $u^\epsilon$.

Next, we want to show the $C^{1,1}$  estimates of $u^{\epsilon}$ are independent of $\epsilon$. Firstly, by maximum principal, $u^{\epsilon_1}\ge u^{\epsilon_2}$ for any $\epsilon_1\le \epsilon_2$. Thus $u^{0}:=\lim\limits_{\epsilon\rightarrow\infty}u^{\epsilon}$ exists. If we could prove uniform $C^{1,1}$ estimates, then   $u^{0}$ is the $C^{1,1}$ solution of  equation \eqref{EquaRing}.
\begin{theorem}
	Let $u^{\epsilon}$ be the smooth $k$-subharmonic solution of \eqref{ApprEquaRing3}. Then there exists a uniform constant $C$ independent of $\varepsilon$ such that
	\begin{align*}
		|u^{\varepsilon}|_{C^{1,1}(\overline \Omega)}\le& C.
	\end{align*}
\end{theorem}
\textbf{In the following subsections, for simplicity, we use $u$ instead of $u^{\epsilon}$.}
\subsection{$C^1$-estimates}
\begin{lemma}
There exists a uniform constant $C$ such that
\begin{align}
	|u|_{C^1(\overline U)}\le C.
	\end{align}
\end{lemma}
\begin{proof}
Let $h$ be the unique solution of the problem
\begin{align}\label{mainequation}
\begin{cases}
\Delta u=0&\quad\text{in }\Omega,\\
h=0&\quad\text{on }\partial \Omega_1,\\
h=1&\quad\text{on }\partial\Omega_0.
\end{cases}
\end{align}
By the maximal principle, we have $\underline{u}\le u\le h$. This gives the uniform $C^{0}$ estimates.

Let $F^{ij}:=\frac{\partial }{\partial u_{ij}}\log H_k(u)=\frac{\partial }{\partial u_{ij}}S_k(\partial\bar\partial u)$.
$$D_\xi=\sum\limits_{i=1}^n(a_i\frac{\partial }{\partial x_i}+b_i\frac{\partial }{\partial y_i})\quad\text{ with }\quad \sum\limits_{i=1}^na^2_i+b_i^2=1.$$ 
Then $$F^{ij}(D_\xi u)_{ij}=0.$$
Thus we have
$$\max\limits_{\overline U}|Du|
=\max\limits_{\partial U}|Du|.$$
Since $\underline{u}\le u^{\varepsilon}\le h$ in $\Omega$ and $\underline{u}= u^{\varepsilon}= h$ on $\p \Omega$, we have
$$h_\nu\leq u^\varepsilon_\nu\leq \underline u_\nu$$
where $\nu$ is the unit outer normal to $\partial \Omega$ (unit inner normal to $\p\Omega_1$ and unit outer normal to $\partial\Omega_0$).
Thus we have
\begin{align}
	\max\limits_{\overline \Omega}|Du|
=\max\limits_{\partial \Omega}|Du|\le C.
\end{align}
\end{proof}
\subsection{Second order estimates}\label{subsection3.2}
\begin{lemma}
	There exists a uniform constant $C$ such that
	\begin{align}
		\max_{\overline U}|D^2 u|\le C.
	\end{align}
\end{lemma}
\begin{proof}
Denote by $D_\xi u=u_\xi$.
Then
$$ L(u_{\xi\xi})=-\frac{\partial^2 S_k^{\frac1k}(\partial\bar\partial u)}{\partial u_{j\bar k}\partial u_{l\bar m}}u_{j\bar k\xi}u_{l\bar m\xi}\geq0.$$
Hence
$$ u_{\xi\xi}(z)\leq \sup_{\partial\Omega}|D^2u|.$$
This implies $ \forall\ i,j=1,\cdots,n,$
\begin{equation*}\begin{aligned}
&u_{x_ix_i},u_{y_iy_i}\leq \sup_{\partial\Omega}|D^2u|,\\
&u_{x_i\pm x_j},u_{x_i\pm y_j},u_{y_i\pm y_j}\leq 2\sup_{\partial\Omega}|D^2u|.
\end{aligned}\end{equation*}
On the other hand, $ \Delta u(z)>0$ implies
$$ u_{x_ix_i},u_{y_iy_i}\geq -(2n-1)(\sup_{\partial\Omega}|D^2u|).$$
Then
\begin{equation*}\begin{aligned}
&\pm u_{x_ix_j}=u_{x_i\pm x_j}-u_{x_ix_i}-u_{x_jx_j}\leq (4n-1)(\sup_{\partial\Omega}|D^2u|),\\
&\pm u_{x_iy_j}=u_{x_i\pm y_j}-u_{x_ix_i}-u_{y_jy_j}\leq (4n-1)(\sup_{\partial\Omega}|D^2u|),\\
&\pm u_{y_iy_j}=u_{y_i\pm y_j}-u_{y_iy_i}-u_{y_jy_j}\leq (4n-1)(\sup_{\partial\Omega}|D^2u|).
\end{aligned}\end{equation*}
Thus we have 
$$\max\limits_{\overline\Omega}|D^2u|\le C_n\max\limits_{\partial\Omega}|D^2u|.$$
So we need to prove the second order estimate on the boundary $\partial \Omega$.  Here we use the method by B. Guan \cite{guan1994cpde,Guanbo2014},  P. F. Guan\cite{gpf2002am} and S. Y. Li \cite{Lisongying2004} .

\emph{Tangential derivative estimates on $\partial \Omega$}\\

Consider a point $p\in \partial \Omega$. Without loss of generality, let $ p $ be the origin. Choose the coordinate $z_1,\cdots,z_n$ such that the $x_n$ axis is the inner normal direction to $\partial\Omega$ at $0$. 
Suppose 
$$ t_1=y_1,\ t_2=y_2,\ \cdots, \ t_n=y_n,\ t_{n+1}=x_1,\ t_{n+2}=x_2, \ \cdots ,\ t_{2n}=x_n.$$
Denote by $t'=(t_1,\cdots,t_{2n-1})$. Then around the origin, $\partial\Omega$ can be represented as a graph
$$ t_n=x_n=\rho(t')=B_{\alpha\beta}t_\alpha t_\beta+O(|t'|^3).$$
Since
$$ (u-\underline u)(t',\rho(t'))=0,$$
we have
$$ (u-\underline u)_{t_\alpha t_\beta}(0)=-(u-\underline u)_{t_n}(0)B_{\alpha\beta},\quad \alpha,\beta=1,\cdots,2n-1.$$
It follows by gradient estimate that
\begin{equation}\label{puretangentialfordirichlet} 
|u_{t_\alpha t_\beta}(0)|\leq C,\quad \alpha,\beta=1,\cdots,2n-1.
\end{equation}

\emph{Tangential-normal derivative estimates on $\partial \Omega$.}\\
 We use Guan's method \cite{guan1994cpde,Guan1999,Guanbo2014} . Our barrier function here is simpler than  before since $u$ is constant on the boundary and the right hand side of the approximating equation is a sufficiently small constant $\epsilon$.

To estimate $ u_{t_\alpha t_n}(0) $ for $\alpha=1,\cdots,2n-1$ and $ u_{t_nt_n}(0)$, we consider the auxiliary function 
$$ v=u-\underline u+td-\frac N2 d^2 $$
on $\Omega_\delta=\Omega\cap B_\delta(0)$ with constant $N,t,\delta$ to be determined later. The following lemma proven in \cite{Guan1999} are needed.

\begin{lemma}\label{lemma4.5-1}
For $N$ sufficiently large and $t,\delta$ sufficiently small, there holds
\begin{equation*}\begin{cases} 
Lv\leq -\frac{\epsilon}4(1+\mathcal F)&\quad\text{in }\Omega_\delta,\\
v\geq 0&\quad\text{on }\partial\Omega,
\end{cases}\end{equation*}
where $\epsilon>0$ is a uniform constant depending only on subsolution $\underline u$ restricted in a small neighborhood of $\partial\Omega$.
\end{lemma}

The following three lemmas was proven by Guan in \cite{Guanbo2014}.
\begin{lemma}\label{lemma719::1}
Let $F^{i\bar j}=\frac{\partial}{\partial u_{i\bar j}}S_k^\frac1k(\partial\bar\partial u) $. Then there is an index $r$ such that 
\begin{align}\label{eq719::4}\sum_{l=1}^{n-1}F^{i\bar j}u_{i\bar l}u_{l\bar j}\geq \frac12 \sum_{i\neq r}S_k^{\frac1k-1}(\lambda)S_{k-1}(\lambda|i)\lambda_i^2,\end{align} where $\lambda=(\lambda_1,\cdots,\lambda_n)$ are the eigenvalues of $u_{i\bar j}$.
\end{lemma}
\begin{lemma}\label{lemma719::2}
Suppose $\lambda\in\overline\Gamma_k$. If $\lambda_r<0$, then
$$\sum_{i\neq r}S_{k}^{\frac1k-1}(\lambda)S_{k-1}(\lambda|i)\lambda_i^2\geq \frac1n \sum_{i=1}^nS_{k}^{\frac1k-1}(\lambda)S_{k-1}(\lambda|i)\lambda_i^2.$$
\end{lemma}
\begin{lemma}\label{lemma719::3}
Suppose $\lambda\in\overline\Gamma_k$. Then for any $r=1,\cdots,n$ and $\varepsilon>0$, 
\begin{align}\label{eq719::6}\sum_{i=1}^nS_k^{\frac 1k-1}(\lambda)S_{k-1}(\lambda|i)|\lambda_i|\leq \varepsilon\sum_{i\neq r}S_k^{\frac1k-1}(\lambda)S_{k-1}(\lambda|i)\lambda_i^2+\frac{C}{\varepsilon}\sum_{i=1}^nS_k^{\frac1k-1}(\lambda)S_{k-1}(\lambda|i)+Q(r),\end{align}
where $Q(r)=S_k^\frac1k(\lambda)-(C_n^k)^\frac1k$ if $\lambda_r\geq 0$ and $Q(r)=0$ if $\lambda_r<0$.
\end{lemma}

At any boundary point $p\in\partial\Omega$, we may choose coordinates $z_1,\cdots,z_n$ with the origin $p$ such that the positive $x_n$ axis is the interior normal direction to $\partial\Omega$ at $p$. Let $\varrho$ be a defining function of $\Omega$, that is $\varrho<0$ in $\Omega$, $\varrho=0$, $D\nu\varrho=1$, on $\partial\Omega$, where $\nu$ is a unit outer normal to $\partial\Omega$. 
We may assume that, $\frac{\partial\varrho}{\partial x_j}(0)=0$ for $1\leq i\leq n-1$ and $\frac{\partial \varrho}{\partial y_j}(0)=0$ for all $1\leq i\leq n$. Moreover, around the orgin, we can write
$$\varrho(z)=-x_n+\mathrm{Re}\sum_{i,j=1}^n\varrho_{ij}(0)z_iz_j+\sum_{i,j=1}^n\varrho_{i\bar j}(0)z_i\bar z_j+Q(z),$$
where $|Q(z)|\leq C|z|^3$.
Let
$$ t_i=y_i,\quad\ i=1,\cdots, n,\quad t_{n+i}=x_i,\quad\ i=1,\cdots,n.$$
Let$$a_\alpha(z)=-\frac{\frac{\partial \varrho}{\partial t_\alpha}}{\frac{\partial\varrho}{\partial x_n}},\quad 1\leq \alpha\leq 2n-1.$$
Then 
$$a_\alpha(0)=0.$$
So
$ T=\frac{\partial}{\partial t_\alpha}+a_\alpha\frac{\partial }{\partial x_n}$ is a tangential vector to $\partial\Omega$ near the origin. We write
$$ a_\alpha(z)=\sum_{\beta=1}^{2n-1}b_{\alpha\beta}t_\beta+b_{\alpha}x_n+O(|t|^2+x_n^2),\quad z\in \overline\Omega \text{ near } 0.$$
And let
$$ T_\alpha=\frac{\partial }{\partial t_\alpha}+\sum_{\beta=1}^{2n-1}b_{\alpha\beta}t_\beta\frac{\partial }{\partial x_n}.$$
Then
$$T=T_\alpha+b_\alpha x_n\frac{\partial }{\partial x_n}+O(|z|^2)\frac{\partial }{\partial x_n}.$$
So 
$$T_\alpha(u-\underline u)=O(|t|^2),\quad\text{on }\partial \Omega.$$
Note that 
$$ \partial _it_{\beta}=\begin{cases}-\frac{\sqrt{-1}}2\delta_{i\beta},\quad 1\leq \beta\leq n,\\\frac12\delta_{i\beta-n},\quad \beta>n.\end{cases}$$
and
$$ \partial _{\bar j}t_{\beta}=\begin{cases}\frac{\sqrt{-1}}2\delta_{j\beta},\quad 1\leq \beta\leq n,\\\frac12\delta_{j\beta-n},\quad \beta>n.\end{cases}$$
We then have 
\begin{equation*}\begin{aligned}
LT_{\alpha}(u-\underline u):=&T_{\alpha}f-LT_{\alpha}\underline u+\sum_{\beta=1}^{2n-1}b_{\alpha\beta}F^{i\bar j}(t_{\beta,i}u_{x_n\bar j}+t_{\beta,\bar j}u_{x_n i})\\
=&T_{\alpha}f-LT_{\alpha}\underline u+2\sum_{\beta=1}^{2n-1}b_{\alpha\beta}F^{i\bar j}(t_{\beta,i}u_{n\bar j}+t_{\beta,\bar j}u_{\bar n i})+\sqrt{-1}\sum_{\beta=1}^{2n-1}b_{\alpha\beta}F^{i\bar j}(t_{\beta,i}u_{y_n\bar j}-t_{\beta,\bar j}u_{y_n i})\\
\geq&-C\Big(1+\sum_{i=1}^n\frac{S_{k-1}(\lambda|i)}{S_k(\lambda)}+\frac{S_{k-1}(\lambda|i)|\lambda_i|}{S_k(\lambda)}\Big)-\frac14F^{i\bar j}(u_{y_ni}-\underline u_{y_ni})(u_{y_n\bar j}-\underline u_{y_n\bar j}),
\end{aligned}\end{equation*}
and
\begin{align*}
	&L\big((u_{y_n}-\underline u_{y_n})^2+\sum_{l=1}^{n-1}|u_l-\underline u_l|^2)\\=&2F^{i\bar j}(u_{y_ni}-\underline u_{y_ni})(u_{y_n\bar j}-\underline u_{y_n\bar j})+\sum_{l=1}^{n-1}F^{i\bar j}\big((u_{li}-\underline u_{li})(u_{\bar l\bar j}-\underline u_{\bar l\bar j})+(u_{l\bar j}-\underline u_{l\bar j})(u_{\bar li}-\underline u_{\bar li})\big)\\
	&+2(u_{y_n}-\underline u_{y_n})F^{i\bar j}(u_{y_ni\bar j}-\underline u_{y_ni\bar j})+\sum_{l=1}^{n-1}\big((u_l-\underline u_l)F^{i\bar j}(u_{\bar li\bar j}-\underline u_{\bar li\bar j})+(u_{\bar l}-\underline u_{\bar l})F^{i\bar j}(u_{li\bar j}-\underline u_{li\bar j})\big)\\
	\geq&2F^{i\bar j}(u_{y_ni}-\underline u_{y_ni})(u_{y_n\bar j}-\underline u_{y_n\bar j})+\sum_{l=1}^{n-1}F^{i\bar j}u_{l\bar j}u_{\bar li}-C\Big(1+\sum_{i=1}^n\frac{S_{k-1}(\lambda|i)}{S_k(\lambda)}+\frac{S_{k-1}(\lambda|i)|\lambda_i|}{S_k(\lambda)}\Big).
\end{align*}

Let
$$\Psi=A_1v+A_2|z|^2-A_3\big((u_{y_n}-\underline u_{y_n})^2+\sum_{l=1}^{n-1}|u_l-\underline u_l|^2).$$
By Lemma \ref{lemma4.5-1}, Lemma \ref{lemma719::1} and Lemma \ref{lemma719::3}, we see that 
\begin{align*}
L(\Psi\pm T_\alpha (u-\underline u))\leq 0\quad\text{in }\Omega_\delta
\end{align*}
and 
\begin{align*}
\Psi\pm T_\alpha(u-\underline u)\geq 0\quad\text{on }\partial\Omega_\delta,
\end{align*}
when $A_1\gg A_2\gg A_3\gg1$.
Therefore
$$ |u_{t_\alpha x_n}|\leq C.$$
In particular, from \eqref{puretangentialfordirichlet}, we know
$$|u_{y_ny_n}|\leq C.$$

\emph{Double normal derivative estimates on $\partial \Omega$ }\\
	For any fixed $p\in\partial \Omega$, we choose the coordinate such that $p=0$, $\partial \Omega\bigcap B_r(0)=(t', \varphi(t'))$ and $\nabla \varphi(0)=0$.
	
\emph{\textbf{Case 1:} $x_0\in\partial \Omega_0$}
Let $\rho_0$ be a defining function of $\Omega_0$ which is strictly plurisubharmonic in a neighborhood of $\Omega_0$. So
$$\rho_0(t',\varphi(t'))=0\quad\text{on }\partial\Omega_0.$$
Then we have
$$\rho_{0,t_\alpha t_\beta}(0)=-\rho_{0,t_{2n}}(0)\varphi_{t\alpha t_\beta}(0)\quad 1\leq \alpha,\beta\leq 2n-1.$$
On the other hand, we have
\begin{align*}
	u_{\alpha\beta}(0)=-u_{t_{2n}}(0)\varphi_{\alpha\beta}(0)\quad 1\leq \alpha,\beta\leq 2n-1.
	\end{align*}
Thus
$$u_{t_\alpha t_\beta}(0)=\frac{u_{t_{2n}}(0)}{\rho_{0,t_{2n}}(0)}\rho_{0,t\alpha t_\beta}(0)\quad 1\leq \alpha,\beta\leq 2n-1.$$
and
$$u_{i\bar j}(0)=\frac{u_{t_{2n}}(0)}{\rho_{0,t_{2n}}(0)}\rho_{0,i\bar j}(0)\geq c\rho_{0,i\bar j}(0)>0.$$
Since $\rho_0$ is strictly plurisubharmonic in $\overline\Omega_0$, we have
\begin{align}\label{0717c1}
	S_{k-1}(\{u_{i\bar j}(0)\}_{1\leq i,j\leq n-1})\ge c^{k-1}S_{k-1}(\{\rho_{0,i\bar j}(0)\}_{1\leq i,j\leq n-1})\ge c_1>0.
	\end{align}

\emph{\textbf{Case 2:} $x_0\in\partial \Omega_1$}

Note that $u\geq \underline u$ near $\partial\Omega_1$, $u=\underline u$ and $0<\underline u_\nu\leq u_\nu$ on $\partial\Omega_1$, $\nu$ is the unit outer normal to $\partial\Omega_1$, there exists a smooth function $g$ such that $u=g\underline u$ near $\partial\Omega_1$, and $g\geq 1$ outside of $\Omega$ nearby $\partial\Omega_1$. So $\forall\,1\leq i,j\leq n-1$,
\begin{align*}
u_{i\bar j}(0)=g_{i\bar j}(0)\underline u(0)+g_i(0)\underline u_{\bar j}(0)+g_{\bar j}(0)\underline u_i(0)+g(0)\underline u_{i\bar j}(0).
\end{align*}
Note that $\underline u=\tau\rho_1$ near $\partial\Omega_1$, where $\rho_1$ is a given strictly plurisubharmonic function in a neighborhood $\Omega$, $\tau$ the a constant independent of $\varepsilon$ and $R$ as taken in Lemma \ref{subu0720}. We also have
\begin{align}\label{0717c2}
S_{k-1}(\{u_{i\bar j}(0)\}_{1\leq i,j\leq n-1})=&\tau^{k-1}g^{k-1}(0)S_{k-1}(\{\rho_{1,i\bar j}(0)\}_{1\leq i,j\leq n-1})\nonumber\\
\geq& \tau^{k-1}g_0^{k-1}C_n^{k-1}(C_n^k)^{\frac{1-k}k}\min_{\partial\Omega}S_k^\frac{k-1}k(\partial\bar \partial \rho_1):=c_1>0.
\end{align}

 Let $c_0=\min\{c_1, c_2\}$ (see \eqref{0717c1} and \eqref{0717c2}),  we have
\begin{align*}
u_{n\bar n}(0)c_0\le& u_{n\bar n}(0)S_{k-1}(\{u_{i\bar j}(0)\}_{1\leq i,j\leq n-1})\\
=&S_k(\{u_{i\bar j}(0)\}_{1\leq i,j\leq n})-S_{k}(\{u_{i\bar j}(0)\}_{1\leq i,j\leq n-1})
+\sum_{i=1}^{n-1}|u_{i\bar n}|^2S_{k-2}(\{u_{i\bar j}(0)\}_{1\leq i,j\leq n-1})\\
\le& C.
\end{align*}
Then we obtain
\begin{align*}
	u_{n\bar n}(0)\le C,
	\end{align*}
where $C$ is a uniform constant.
  On the other hand, $u_{n\bar n}(0)\ge \sum\limits_{i=1}^{n-1}u_{\alpha\bar \alpha}(0)\ge -C$. In conclusion, we have $|u_{n\bar n}(0)|\le C$.

  In conclusion, we get the uniform $C^2$ estimate.
\end{proof}

\subsection{Proof of Theorem \ref{ring0720}}

The uniqueness follows from the comparison principle for $k$-subharmonic solutions of complex $k$-Hessian equations in Lemma \ref{comparison0718} by Blocki \cite{BlockiAIF}.

For the existence part,  since $u^{\epsilon}$ is increasing on $\epsilon$,  $u^{0}:=\lim\limits_{\epsilon\rightarrow 0}u^{\epsilon}$ exits. Since  $|u^{\epsilon}|_{C^2(\overline \Omega)}\le C$, there exists a subsequence $u^{\epsilon_i}$  converges to $u^0 $ in $C^{1,\alpha}$ on $\overline \Omega$ and $u^0\in C^{1,1}(\overline \Omega)$.

\section{Solving the approximating equation in $\Sigma_R:=B_R\setminus \Omega$.}\label{section4}
We always assume $\Omega$ is a smooth, strongly  pseudoconvex domain containing the orgin and $\Omega$ is holomorphically convex in a ball. Recall that we always assume $B_{r_0}\subset\subset\Omega\subset\subset B_{R_0}\subset\subset B_{S_0}$ and $\Omega$ is holomorphically convex in $B_{S_0}$. 

Since the Green function in this case is $-|z|^{2-\frac{2n}k}$, we want to solve the following complex $k$-Hessian equation .
\begin{equation}\label{case1Equa}
\begin{cases}
(dd^cu)^k\wedge\omega^{n-k}=0&\quad\text{in }\Omega^c:=\mathbb C^n\backslash\overline\Omega,\\
u=-1&\quad\text{on }\partial\Omega,\\
u(z)\rightarrow 0&\quad\text{as }|z|\rightarrow\infty.
\end{cases}
\end{equation}
By scaling of $z$, we consider \eqref{case1Equa} with $B_t\subset\subset\Omega\subset\subset B_1\subset\subset B_{1+s}$, where $t=\frac{r_0}{R_0}$, $s=\frac{S_0}{R_0}-1$.

\subsection{Construction of the approximating equation. }
Let $w^\varepsilon$ be a approximation of the Green function $-|z|^{2-\frac{2n}k}$
$$w^\varepsilon(z)=-\bigg(\frac{|z|^2+\varepsilon^2}{1+\varepsilon^2}\bigg)^{1-\frac nk}.$$
We have
$$f^\varepsilon:=H_k(w^\varepsilon)=S_k(w^\varepsilon_{i\bar j})=C_n^k(\frac nk-1)^k\varepsilon^2(1+\varepsilon^2)^{n-k}(|z|^2+\varepsilon^2)^{-n-1}.$$

It is clear that $\rho_0=|z|^2-(1+s)^2$ is a plurisubharmonic defining function of $B_{1+s}$. Let $\rho_1$ be a defining function of $\Omega$ such that $\rho_1$ is plurisubharmonic in a neighborhood $U$ of $\Omega$.

By Lemma \ref{subu0720}, there is a smooth plurisubharmonic function $\rho$ solving
\begin{equation}\label{underlineu071911-2}
\begin{cases}
H_k(\rho)\geq \epsilon_0,&\quad\text{in }B_{1+s}\backslash\overline\Omega,\\
\rho=\mathrm\tau\rho_1,&\quad\text{near }\partial\Omega,\\
\rho=1+\mathrm K\rho_0,&\quad\text{near }\partial B_{1+s},\\
\end{cases}
\end{equation}
Let $\varphi=(1-(1+\frac{s^2}{16+s^2})^{1-\frac{n}k})\rho-1$. In $B_{1+s}\backslash B_{1+\frac s2}$, $\forall\ \varepsilon\leq \varepsilon_0$, $\varepsilon_0<\frac{s^2}8$, 
$$w^\varepsilon\geq -\bigg(\frac{(1+\frac s2)^2}{1+\varepsilon_0^2}\bigg)^{1-\frac nk}>-\bigg(1+\frac{s^2}{8+s^2}\bigg)^{1-\frac nk}. $$
So 
$$w^\varepsilon-\varphi>\bigg(1+\frac{s^2}{16+s^2}\bigg)^{1-\frac nk}-\bigg(1+\frac{s^2}{8+s^2}\bigg)^{1-\frac nk}\quad\text{in }B_{1+s}\backslash B_{1+\frac s2}.$$
Let $V$ be a neighborhood of $\Omega$, $\Omega\subset\subset V$, then 
$$w^\varepsilon\leq -1\quad\text{and}\quad \varphi\geq (1-(1+\frac{s^2}{16+s^2})^{1-\frac{n}k})\inf_{B_1\backslash V}\rho-1\quad\text{in }B_1\backslash V.$$
So
$$w^\varepsilon-\varphi\leq (1-(1+\frac{s^2}{16+s^2})^{1-\frac{n}k})\inf_{B_1\backslash V}\rho \quad\text{in }B_1\backslash V.$$

Apply Lemma \ref{Guan2002} with $w^\varepsilon$, $\varphi$ and $\delta<\min\bigg\{\bigg(1+\frac{s^2}{16+s^2}\bigg)^{1-\frac nk}-\bigg(1+\frac{s^2}{8+s^2}\bigg)^{1-\frac nk},(1-(1+\frac{s^2}{16+s^2})^{1-\frac{n}k})\inf\limits_{B_1\backslash V}\rho\bigg\}$, we obtain a smooth $k$-subharmonic function $\underline u^\varepsilon$ such that $\underline u^\varepsilon=w^\varepsilon$ in $\mathbb C^n\backslash B_{1+\frac s2}$, $\underline u^\varepsilon=\varphi$ in $B_1\backslash \Omega$ and $\underline u^\varepsilon\geq \max\{\varphi,w^\varepsilon\}$ in $\Omega^c$. Moreover, by the convacity of $S_k^\frac1k$,
$$H_k^\frac1k(\underline u^\varepsilon)\geq \frac{1+t(z)}2H_k^\frac1k(\varphi)+\frac{1-t(z)}2H_k^\frac1k(w^\varepsilon)\quad\text{in }\{|\varphi-w^\varepsilon|<\delta\}.$$
If we take $\varepsilon_0<\min\{1,2^{k-n}t^{-2(n+1)}(C_n^k)^{-1}(\frac nk-1)^k(1-(1+\frac{s^2}{16+s^2})^{1-\frac{n}k})\epsilon_0\}$, then for any $\varepsilon\leq \varepsilon_0$, $f^\varepsilon<\epsilon_0$. So we obtain 
$$H_k(\underline u^\varepsilon)\geq f^\varepsilon\quad\text{in }\Omega^c.$$

In conclusion, for sufficient small $\varepsilon$, we can construct a smooth, strictly $k$-subharmonic function $\underline u^\varepsilon$ as follows
\begin{lemma}
For any $\varepsilon\in(0,\varepsilon_0)$, $\varepsilon_0<\frac{s^2}8$, there exists a strictly $k$-subharmonic function $\underline u^\varepsilon\in C^\infty(\mathbb C^n\backslash \Omega)$ satisfying
\begin{align*}
\underline u^\varepsilon=
\begin{cases}
w^\varepsilon&\quad\text{in }\mathbb C^n\backslash B_{1+\frac s2},\\
(1-(1+\frac{s^2}{16+s^2})^{1-\frac{n}k})\rho-1&\quad\text{in }B_1\backslash \Omega,
\end{cases}
\end{align*}
$$\underline u^\varepsilon\geq \max\{w^\varepsilon,(1-(1+\frac{s^2}{16+s^2})^{1-\frac{n}k})\rho-1 \}\quad\text{in }B_{1+\frac s2}\backslash B_{1},$$
and 
$$H_k(\underline u^\varepsilon)\geq f^\varepsilon\quad\text{in }\Omega^c,$$
where $\rho$ is a function satisfying \eqref{underlineu071911-2}.
\end{lemma}

By the preliminaries above in this section, we are able to construct the approximation equations for $\varepsilon\in (0,\varepsilon_0)$ and $R>1+s$.
\begin{equation}\label{approxeq2}
\begin{cases}
H_k(u^{\varepsilon,R})=f^\varepsilon&\quad\text{in }\Sigma_R:=B_R\backslash\Omega,\\
u^{\varepsilon,R}=\underline u^\varepsilon&\quad\text{on }\partial\Sigma_R.
\end{cases}
\end{equation}
Since $\underline u^\varepsilon$ is a subsolution,  by Li \cite{Lisongying2004}, \eqref{approxeq2} has a strictly $k$-subharmonic solution $u^{\varepsilon,R}\in C^\infty (\overline \Sigma_R)$. Our goal is to establish uniform $C^2$ estimates of $u^{\varepsilon,R}$, which is independent of $\varepsilon$ and $R$. We prove the following 
\begin{theorem}\label{kindofest}
For sufficient small $\varepsilon$ and sufficient large $R$, $u^{\varepsilon,R}$ satisfies 
\begin{align*}
&C^{-1}|z|^{2-\frac{2n}k}\leq -u^{\varepsilon,R}(z)\leq C|z|^{2-\frac{2n}k},\\
&|Du^{\varepsilon,R}(z)|\leq C|z|^{1-\frac{2n}k},\\
&|\p\bar\p  u^{\varepsilon,R}(z)|\leq C|z|^{-\frac{2n}k},\\
&|D^2u^{\varepsilon,R}(z)|\leq C.
\end{align*}
where $C$ is a uniform constant which is independent of $\varepsilon$ and $R$.
\end{theorem}
In the next subsections, we will prove uniform $C^2$-estimates of solutions to equation \eqref{approxeq2}. The key point is that these estimates are independent of $\varepsilon$ and $R$.

\subsection{$C^0$ estimates.} 
Since $\underline u^\varepsilon$ is a subsolution to \eqref{approxeq2}, we obtain that
$$u^{\varepsilon,R}\geq \underline u^\varepsilon\geq -\Big(\frac{|z|^2+\varepsilon^2}{1+\varepsilon^2}\Big)^{1-\frac nk}\geq -(1+\varepsilon_0^2)^{\frac nk-1}|z|^{2-\frac {2n}k}.$$
For any $R'\geq R\geq 1+s$, let $u^{\varepsilon,R}$ and $u^{\varepsilon,R'}$ be solutions to \eqref{approxeq2} on $\Sigma_R$ and $\Sigma_{R'}$ respectively. We have
$$u^{\varepsilon,R}=\underline u^\varepsilon\leq u^{\varepsilon,R'}\quad\text{on }\partial B_R.$$
By Lemma \ref{comparison0718}, 
$$ u^{\varepsilon,R}\leq u^{\varepsilon,R'},\quad\text{in }\Sigma_R.$$
On the other hand, choose 
$R_1:=\max\bigg\{1+s,\frac{t\varepsilon_0}{\sqrt{1-t^2}}\bigg\}$. Then for any $R\geq R_1$,
\begin{equation*}\begin{cases}
H_k(-t^{\frac{2n}k-2}|z|^{2-\frac{2n}k})=0<f^\varepsilon=H_k(u^{\varepsilon,R})&\quad\text{in }\Sigma_R,\\
u^{\varepsilon,R}=-1\leq -t^{\frac{2n}k-2}|z|^{2-\frac{2n}k}&\quad\text{on }\partial\Omega,\\
u^{\varepsilon,R}=-\big(\frac{R^2+\varepsilon^2}{1+\varepsilon^2}\big)^{1-\frac nk}\leq -t^{\frac{2n}k-2}R^{2-\frac{2n}k}&\quad\text{on }\partial B_R,
\end{cases}\end{equation*}
Using Lemma \ref{comparison0718} again, we have
$$ u^{\varepsilon,R}\leq -t^{\frac{2n}k-2}|z|^{2-\frac{2n}k}\quad\text{in }\Sigma_R.$$
So we have, for any $R'>R\geq R_1$, 
$$-(1+\varepsilon_0^2)^{\frac nk-1}|z|^{2-\frac {2n}k}\leq u^{\varepsilon,R}(z)\leq u^{\varepsilon,R'}(z)\leq -t^{\frac{2n}k-2}|z|^{2-\frac{2n}k},\quad z\in \Sigma_R.$$

\subsection{Gradient estimates}
In this subsection, we prove the global gradient estimate. The key point is that the estimate here does not depend on $\varepsilon$ and $R$. We also prove that the positive lower bound of the  gradient of the solution. 
\subsubsection\textbf{\emph{{Reducing global gradient estimates to boundary gradient estimates.}}}
This part is the key part of gradient estimates. The point in here is that the gradient estimate is independent of the approximating process.  This estimates is motivated by B. Guan \cite{gb2007imrn}.
\begin{theorem}\label{Gradientest}
Let $u$ be the solution of the approximating equation \eqref{approxeq2}.
Denote by
\begin{align}\label{threecase}
P=|Du|^2(-u)^{-\frac{2n-k}{n-k}}.
\end{align}
then we have the following gradient estimate
\begin{align}\label{Gradientestimate}
\max_{\Sigma_R}P\leq \max\bigg\{\max_{\partial \Sigma_R}P,\Big(\frac{2(n-k)}{k(2n-k)}\Big)^2(-u)^{-\frac{k}{n-k}}|D\log f^\varepsilon|^2\bigg\}.
\end{align}

\end{theorem}
\begin{proof}For simplicity, we use $f$ instead of $f^\varepsilon$ during the proof.

Let $a=\frac{2n-k}{n-k}$.
Select the auxiliary function 
$$\varphi=\log P=\log |Du|^2-a\log(-u).$$
Suppose $\varphi$ obtain its maximum at $z_0\in \Sigma_R$. We can choose the holomorphic coordinate such that $\{u_{i\bar j}\}(z_0)$ is diagonal. Denote by $\lambda_i=u_{i\bar i}(z_0)$. The following computations are at $z_0$.
\begin{equation*}
0=\varphi_i=\frac{|Du|^2_i}{|Du|^2}-a\frac{u_i}{u}=\frac{u_lu_{\bar li}+u_{li}u_{\bar l}}{|Du|^2}-a\frac{u_i}u=\frac{u_i\lambda_i+u_{li}u_{\bar l}}{|Du|^2}-a\frac{u_i}u.
\end{equation*}
Then we have the observation
\begin{align}\label{key::1}
a\frac{|u_i|^2}{u}=\frac{|u_i|^2\lambda_i}{|Du|^2}+\sum\limits_{l=1}^n\frac{u_{li}u_{\bar l}u_{\bar i}}{|Du|^2},\quad\forall\ i=1,\cdots,n,
\end{align}
which implies $\sum_{l=1}^nu_{li}u_{\bar l}u_{\bar i}$ is real at $z_0$.
Denote by $F^{ij}=\frac{\partial}{\partial u_{i\bar j}}S_k(\partial\bar\partial u)$. By direct computation, we can  get
\begin{equation*}\begin{aligned}
0\geq & F^{i\bar j}\varphi_{i\bar j}=F^{i\bar j}\cdot\bigg(\frac{|Du|^2_{i\bar j}}{|Du|^2}-\frac{|Du|^2_i|Du|^2_{\bar j}}{|Du|^4}-a\frac{u_{i\bar j}}{u}+a\frac{u_iu_{\bar j}}{u^2}\bigg)\\
=&F^{i\bar j}\cdot\bigg(\frac{|Du|^2_{i\bar j}}{|Du|^2}-(1-\frac1a)\frac{|Du|^2_i|Du|^2_{\bar j}}{|Du|^4}-a\frac{u_{i\bar j}}{u}\bigg)\\
=&\frac{2\mathrm{Re}\{u_lf_{\bar l}\}}{|Du|^2}-akf\frac{|Du|^2}{u}+\sum\limits_{i,l=1}^n\frac{S_{k-1}(\lambda|i)|u_{li}|^2}{|Du|^2}+\sum\limits_{i=1}^n\frac{S_{k-1}(\lambda|i)\lambda_i^2}{|Du|^2}\\
&-\frac{n}{2n-k}\sum_{i=1}^n\frac{|u_i|^2}{|Du|^4}S_{k-1}(\lambda|i)\lambda_i^2-\frac{n}{2n-k}\sum_{i=1}^nS_{k-1}(\lambda|i)\frac{|\sum_{l=1}^nu_{\bar l}u_{li}|^2}{|Du|^4}-\frac{2n}{2n-k}\sum_{i=1}^nS_{k-1}(\lambda|i)\lambda_i\frac{\sum_{l=1}^nu_{li}u_{\bar l}u_{\bar i}}{|Du|^4}.
\end{aligned}\end{equation*}

We claim 
\begin{align}\label{claim::1}	 
\mathcal{E}:=\sum_{i=1}^n\bigg(\sum_{l=1}^nS_{k-1}(\lambda|i)|u_{li}|^2&+S_{k-1}(\lambda|i)\lambda_i^2-\frac{n}{2n-k}\frac{|u_i|^2}{|Du|^2}S_{k-1}(\lambda|i)\lambda_i^2\nonumber\\&-\frac{n}{2n-k}S_{k-1}(\lambda|i)\frac{|\sum_{l=1}^nu_{\bar l}u_{li}|^2}{|Du|^2}-\frac{2n}{2n-k}S_{k-1}(\lambda|i)\lambda_i\frac{\sum_{l=1}^nu_{li}u_{\bar l}u_{\bar i}}{|Du|^2}\bigg)\geq 0.
\end{align}
Then
\begin{equation*}\begin{aligned}
0\geq|Du|^2F^{i\bar j}\varphi_{i\bar j}
\geq 2\mathrm{Re}\{u_lf_{\bar l}\}-akf\frac{|Du|^2}{u}
\geq-2|Du||Df|-akf\frac{|Du|^2}{u}.
\end{aligned}\end{equation*}
It follows that
$$ |Du|\leq \frac{2}{ak}(-u)|D\log f|=\frac{2(n-k)}{k(2n-k)}(-u)|D\log f|. $$
Thus
\begin{align*}
|Du|^2(-u)^{-a}\leq \bigg(\frac{2(n-k)}{k(2n-k)}\bigg)^2(-u)^{2-a}|D\log f|^2.
\end{align*}

Now we prove the Claim \eqref{claim::1}.
Since 
\begin{equation*}\begin{aligned}
\sum_{i=1}^nS_{k-1}(\lambda|i)\lambda_i^2=&S_1f-(k+1)S_{k+1}=\sum_{i=1}^n \frac{|u_i|^2}{|Du|^2}\big(S_1f-(k+1)S_{k+1}\big)\\
=&\sum_{i=1}^nf\frac{|u_i|^2}{|Du|^2}\bigg(\lambda_i+S_{1}(\lambda|i)-(k+1)\frac{S_{k}(\lambda|i)}{S_{k-1}(\lambda|i)}\bigg)\\&+\sum_{i=1}^n\frac{|u_i|^2}{|Du|^2}\bigg((k+1)\frac{S^2_{k}(\lambda|i)}{S_{k-1}(\lambda|i)}-(k+1)S_{k+1}(\lambda|i)\bigg),
\end{aligned}\end{equation*}
we have
\begin{equation*}\begin{aligned}
\mathcal{E}
=&\sum_{i=1}^nf\frac{|u_i|^2}{|Du|^2}\bigg(\lambda_i+S_{1}(\lambda|i)-(k+1)\frac{S_{k}(\lambda|i)}{S_{k-1}(\lambda|i)}\bigg)+\sum_{i=1}^n\frac{|u_i|^2}{|Du|^2}\bigg((k+1)\frac{S^2_{k}(\lambda|i)}{S_{k-1}(\lambda|i)}-(k+1)S_{k+1}(\lambda|i)\bigg)\\
&+\sum_{i,l=1}^nS_{k-1}(\lambda|i)|u_{li}|^2-\frac{n}{2n-k}\sum_{i=1}^n\frac{|u_i|^2}{|Du|^2}S_{k-1}(\lambda|i)\lambda_i^2-\frac{n}{2n-k}S_{k-1}\sum_{i=1}^n(\lambda|i)\frac{|\sum_{l=1}^nu_{\bar l}u_{li}|^2}{|Du|^2}\\
&-\frac{2n}{2n-k}\sum_{i,l=1}^nS_{k-1}(\lambda|i)\lambda_i\frac{u_{li}u_{\bar l}u_{\bar i}}{|Du|^2}\\
:=&\bigg(\sum_{i\in G}+\sum_{i\in H}\bigg)\mathrm T_i,
\end{aligned}\end{equation*}
in which
$$ G=\{i\mid\lambda_i\geq 0\}\quad H=\{i\mid\lambda_i<0\},$$
and
\begin{equation*}\begin{aligned}
\mathrm T_i=&f\frac{|u_i|^2}{|Du|^2}\bigg(\lambda_i+S_{1}(\lambda|i)-(k+1)\frac{S_{k}(\lambda|i)}{S_{k-1}(\lambda|i)}\bigg)+\frac{|u_i|^2}{|Du|^2}\bigg((k+1)\frac{S^2_{k}(\lambda|i)}{S_{k-1}(\lambda|i)}-(k+1)S_{k+1}(\lambda|i)\bigg)\\
&+\sum_{l=1}^nS_{k-1}(\lambda|i)|u_{li}|^2-\frac{n}{2n-k}\frac{|u_i|^2}{|Du|^2}S_{k-1}(\lambda|i)\lambda_i^2-\frac{n}{2n-k}S_{k-1}(\lambda|i)\frac{|\sum_{l=1}^nu_{\bar l}u_{li}|^2}{|Du|^2}\\&-\frac{2n}{2n-k}\sum_{l=1}^nS_{k-1}(\lambda|i)\lambda_i\frac{u_{li}u_{\bar l}u_{\bar i}}{|Du|^2}.
\end{aligned}\end{equation*}
We will prove in the following that $ \forall\  i$, $\mathrm T_i\geq0$.

\noindent\textbf{Case 1. }$i\in H$. Let
$$ \mathrm T_{i}=\mathrm A+\mathrm B,$$
where
\begin{equation*}\begin{aligned}
\mathrm A:=&f\frac{|u_i|^2}{|Du|^2}\bigg(\lambda_i+S_{1}(\lambda|i)-(k+1)\frac{S_{k}(\lambda|i)}{S_{k-1}(\lambda|i)}\bigg)\\
&+\frac{|u_i|^2}{|Du|^2}\bigg((k+1)\frac{S^2_{k}(\lambda|i)}{S_{k-1}(\lambda|i)}-(k+1)S_{k+1}(\lambda|i)\bigg)-\frac{n}{n-k}\frac{|u_i|^2}{|Du|^2}S_{k-1}(\lambda|i)\lambda_i^2,
\end{aligned}\end{equation*}
and
\begin{equation*}\begin{aligned}
 \mathrm B:=\bigg(\frac{n}{n-k}-\frac{n}{2n-k}\bigg)&\frac{|u_i|^2}{|Du|^2}S_{k-1}(\lambda|i)\lambda_i^2+\sum_{l=1}^nS_{k-1}(\lambda|i)|u_{li}|^2\\
&-\frac{n}{2n-k}S_{k-1}(\lambda|i)\frac{|\sum_{l=1}^nu_{\bar l}u_{li}|^2}{|Du|^2}-\frac{2n}{2n-k}\sum_{l=1}^nS_{k-1}(\lambda|i)\lambda_i\frac{u_{li}u_{\bar l}u_{\bar i}}{|Du|^2}.
\end{aligned}\end{equation*}
Since
$$f=S_k(\lambda)=S_{k-1}(\lambda|i)\lambda_i+S_{k}(\lambda|i),$$
we have
$$\lambda_i^2=\frac{f^2}{S^2_{k-1}(\lambda|i)}+\frac{S^2_{k}(\lambda|i)}{S^2_{k-1}(\lambda|i)}-\frac{2fS_{k}(\lambda|i)}{S^2_{k-1}(\lambda|i)}=\frac{f\lambda_i}{S_{k-1}(\lambda|i)}+\frac{S^2_{k}(\lambda|i)}{S^2_{k-1}(\lambda|i)}-\frac{fS_{k}(\lambda|i)}{S^2_{k-1}(\lambda|i)}.$$
Then
$$-\frac{n}{n-k}\frac{|u_i|^2}{|Du|^2}S_{k-1}(\lambda|i)\lambda_i^2=f\frac{|u_i|^2}{|Du|^2}\Big(-\frac{n}{n-k}\lambda_i+\frac{n}{n-k}\frac{S_{k}(\lambda|i)}{S_{k-1}(\lambda|i)}\Big)+\frac{|u_i|^2}{|Du|^2}\Big(-\frac{n}{n-k}\frac{S^2_{k}(\lambda|i)}{S_{k-1}(\lambda)}\Big).$$
By (a) and (b) of Proposition \ref{2ineqs}, we have
\begin{equation*}\begin{aligned}
\mathrm A=&f\frac{|u_i|^2}{|Du|^2}\bigg(\lambda_i+S_{1}(\lambda|i)-(k+1)\frac{S_{k}(\lambda|i)}{S_{k-1}(\lambda|i)}-\frac{n}{n-k}\lambda_i+\frac{n}{n-k}\frac{S_{k}(\lambda|i)}{S_{k-1}(\lambda|i)}\bigg)\\
&+\frac{|u_i|^2}{|Du|^2}\bigg((k+1)\frac{S^2_{k}(\lambda|i)}{S_{k-1}(\lambda|i)}-(k+1)S_{k+1}(\lambda|i)-\frac{n}{n-k}\frac{S^2_{k}(\lambda|i)}{S_{k-1}(\lambda)}\bigg)\\
=&f\frac{|u_i|^2}{|Du|^2}\bigg(-\frac{k}{n-k}\lambda_i+S_{1}(\lambda|i)-(k+1-\frac{n}{n-k})\frac{S_{k}(\lambda|i)}{S_{k-1}(\lambda|i)}\bigg)\\
&+\frac{|u_i|^2}{|Du|^2}\bigg((k+1-\frac{n}{n-k})\frac{S^2_{k}(\lambda|i)}{S_{k-1}(\lambda|i)}-(k+1)S_{k+1}(\lambda|i)\bigg)\\
\geq& f\frac{|u_i|^2}{|Du|^2}\bigg(-\frac k{n-k}\lambda_i+\frac{k}{n-1}S_{1}(\lambda|i)\bigg)\geq 0,
\end{aligned}\end{equation*}
where the last inequality is due to the assumption of \text{Case 1.}
Note that
$$ \sum_{l=1}^nS_{k-1}(\lambda|i)|u_{li}|^2-\frac{n}{2n-k}S_{k-1}(\lambda|i)\frac{|\sum_{l=1}^nu_{\bar l}u_{li}|^2}{|Du|^2}\geq \frac{n-k}{2n-k}\sum_{l=1}^nS_{k-1}(\lambda|i)|u_{li}|^2.$$
And
$$ \frac{2n}{2n-k}S_{k-1}(\lambda|i)\lambda_i\frac{\sum_{l=1}^nu_{li}u_{\bar l}u_{\bar i}}{|Du|^2}\leq \frac1\varepsilon \frac{n^2}{(2n-k)^2}\frac{|u_i|^2}{|Du|^2}S_{k-1}(\lambda|i)\lambda_i^2+\varepsilon S_{k-1}(\lambda|i)\frac{|\sum_{l=1}^nu_{\bar l}u_{li}|^2}{|Du|^2}.$$
Take $\varepsilon=\frac{n-k}{2n-k}$, then $\frac1\varepsilon \frac{n^2}{(2n-k)^2}=\frac{n}{n-k}-\frac{n}{2n-k}$. It follows that 
$\mathrm B\geq 0$.
	
\noindent\textbf{Case 2. } $i\in G$. Then let
$$ \mathrm T_i=\mathrm E+\mathrm F,$$
where
\begin{equation*}\begin{aligned}
\mathrm E:=f\frac{|u_i|^2}{|Du|^2}\bigg(\lambda_i&+S_{1}(\lambda|i)-(k+1)\frac{S_{k}(\lambda|i)}{S_{k-1}(\lambda|i)}\bigg)\\
&+\frac{|u_i|^2}{|Du|^2}\bigg((k+1)\frac{S^2_{k}(\lambda|i)}{S_{k-1}(\lambda|i)}-(k+1)S_{k+1}(\lambda|i)\bigg)-\frac{|u_i|^2}{|Du|^2}S_{k-1}(\lambda|i)\lambda_i^2,
\end{aligned}\end{equation*}
and
\begin{equation*}\begin{aligned}
\mathrm F:=\bigg(1-\frac{n}{2n-k}\bigg)\frac{|u_i|^2}{|Du|^2}&S_{k-1}(\lambda|i)\lambda_i^2+\sum_{l=1}^nS_{k-1}(\lambda|i)|u_{li}|^2\\
&-\frac{n}{2n-k}S_{k-1}(\lambda|i)\frac{|\sum_{l=1}^nu_{\bar l}u_{li}|^2}{|Du|^2}-\frac{2n}{2n-k}\sum_{l=1}^nS_{k-1}(\lambda|i)\lambda_i\frac{u_{li}u_{\bar l}u_{\bar i}}{|Du|^2}.
\end{aligned}\end{equation*}
Since $i\in G$, we have $\lambda_i\geq 0$, it follows from \eqref{key::1} that 
$$ \sum_{l=1}^nu_{li}u_{\bar l}u_{\bar i}<0.$$ 
Then
$$ \mathrm F\geq\frac{n-k}{2n-k}\frac{|u_i|^2}{|Du|^2}S_{k-1}(\lambda|i)\lambda_i^2+\sum_{l=1}^nS_{k-1}(\lambda|i)|u_{li}|^2-\frac{n}{2n-k}S_{k-1}(\lambda|i)\frac{|\sum_{l=1}^nu_{\bar l}u_{li}|^2}{|Du|^2}\geq 0.$$
Using (b) of Proposition \ref{2ineqs}, we obtain
\begin{equation*}\begin{aligned}
\mathrm E=&f\frac{|u_i|^2}{|Du|^2}\left(S_{1}(\lambda|i)-k\frac{S_{k}(\lambda|i)}{S_{k-1}(\lambda|i)}\right)+\frac{|u_i|^2}{|Du|^2}\left(k\frac{S^2_{k}(\lambda|i)}{S_{k-1}(\lambda|i)}-(k+1)S_{k+1}(\lambda|i)\right)\\
\geq& \frac{k-1}{n-1}\frac{|u_i|^2}{|Du|^2}S_{1}(\lambda|i)+\frac{k}{n-k}\frac{|u_i|^2}{|Du|^2}\frac{S^2_{k}(\lambda|i)}{S_{k-1}(\lambda|i)}\geq 0.
\end{aligned}\end{equation*}
Hence we complete the proof of claim \eqref{claim::1}.

\end{proof}
\subsubsection{\textbf{{ Boundary gradient estimates.}}}
We always assume $R>>R_1$. To prove the boundary gradient estimates, we will construct upper barriers on $\partial \Omega$ and $\partial B_R$ respectively.

Let $h_1\in C^\infty(\overline \Sigma_{R_1})$ be the solution of the following equation
\begin{equation*}\begin{cases}
\Delta h_1=0&\quad\text{in }\Sigma_{R_1},\\
h_1=-1&\quad\text{on }\partial\Omega,\\
h_1=-t^{\frac{2n}k-2}|z|^{2-\frac{2n}k}&\quad\text{on }\partial B_{R_1}.
\end{cases}\end{equation*}
$u^{\varepsilon,R}$ is $k$-subharmonic in $\Sigma_R$, thus  is subharmonic in $\Sigma_R$. Note that
$$h_1=u^{\varepsilon,R}=-1\quad\text{on }\partial\Omega \quad\text{and}\quad h_1=-t^{\frac{2n}k-2}R_1^{2-\frac{2n}k}\geq u^{\varepsilon,R}\quad\text{on }\partial B_{R_1}.$$
By comparison theorem for the Laplace equation, we obtain
$$u^{\varepsilon,R}\leq h_1\quad\text{in }\Sigma_{R_1}.$$
Let $\nu$ be the unit outer normal to $\partial\Omega$, then
$$\Big(1-(1+\frac{s^2}{16+s^2})^{1-\frac{n}k}\Big)\rho_\nu= \underline u^\varepsilon_\nu\leq u^{\varepsilon,R}_\nu\leq h_{1,\nu}\leq C(h_1)=C(\Omega,t,R_1)\quad\text{on }\partial\Omega,$$ where $\rho$ is defined in \eqref{underlineu071911-2}.
So there is a constant $C$ independent of $\varepsilon$ and $R$ such that 
$$ |Du^{\varepsilon,R}|\leq C, \quad\text{on }\partial\Omega.$$

Let $h_2\in C^{\infty}(\overline {B_R}\backslash B_{\frac R2})$ be a solution to the following equations,
\begin{equation*}\begin{cases}
\Delta h_2=0&\quad\text{ in }B_R\backslash\overline{B_\frac R2},\\
h_2=\underline u^\varepsilon&\quad\text{ on }\partial B_R,\\
h_2=-(2t)^{\frac{2n}k-2}|z|^{2-\frac{2n}k}&\quad\text{ on }\partial B_\frac R2.
\end{cases}\end{equation*}
For any $C^2$ function $g$, set
$$\tilde g=R^{\frac{2n}k-2}g(R\cdot).$$
Then $\tilde h_2(z)=R^{\frac{2n}k-2}h_2(Rz)$ satisfies 
\begin{equation*}\begin{cases}
\Delta\tilde h_2=0&\quad\text{in }B_1\backslash\overline{B_\frac12},\\
\tilde h_2=\tilde{\underline u}^\varepsilon=-\Big(\frac{1+\frac{\varepsilon^2}{R^2}}{1+\varepsilon^2}\Big)^{1-\frac nk}&\quad\text{on }\partial B_1,\\
\tilde h_2=-(2t)^{\frac{2n}k-2}&\quad\text{on }\partial B_\frac12.
\end{cases}\end{equation*}
Note that
$$\tilde h_2=\tilde u^{\varepsilon,R}\quad\text{on }\partial B_1 \quad\text{and}\quad \tilde h_2\geq \tilde u^{\varepsilon,R}\quad\text{on }\partial B_{\frac12}.$$
By comparison theorem, we obtain
$$u^{\varepsilon,R}\leq \tilde h_2\quad\text{in }B_1\backslash B_{\frac12}.$$
Let $\nu$ be the unit outer normal to $B_1$. Then
$$\tilde h_{2,\nu}\leq \tilde u^{\varepsilon,R}_\nu\leq \underline{\tilde u}^{\varepsilon}_\nu\quad\text{on }\partial B_1.$$
Noted that
$$-2^{1-\frac nk}\geq -\bigg(\frac{1+\frac{\varepsilon^2}{R^2}}{1+\varepsilon^2}\bigg)^{1-\frac nk} \geq-\big(\frac1{1+\varepsilon_0^2}\big)^{1-\frac nk},$$
then $\tilde h_2$ is uniformly bounded on $\partial B_1\backslash\overline{B_\frac12}.$ Since the gradient estimate of harmonic function depends only on the domain and $C^0$ norm of boundary value, there is a positive constant independent of $\varepsilon$ and $R$, such that
$$|\tilde h_{2,\nu}|\leq C,\quad\text{on }\partial B_1.$$
On the other hand, since 
$$\tilde{\underline u}^{\varepsilon}=-\bigg(\frac{|z|^2+\frac{\varepsilon^2}{R^2}}{1+\varepsilon^2}\bigg)^{1-\frac nk},\quad\text{in a neighbourhood of }\partial B_1,$$
we have
$$\tilde{\underline u}^{\varepsilon}_\nu	=(\frac nk-1)\bigg(\frac{1+\frac{\varepsilon^2}{R^2}}{1+\varepsilon^2}\bigg)^{-\frac nk}\frac{\bar z\cdot\nu}{1+\varepsilon^2}\quad\text{on }\partial B_1.$$
Hence
$$|D\tilde u^{\varepsilon,R}|\leq C,\quad\text{on }\partial B_1\quad\text{ independent of }\varepsilon\text{ and }R.$$
So we have the $(\varepsilon,R)$-independent estimate
$$|D u^{\varepsilon,R}|\leq CR^{1-\frac{2n}k},\quad\text{on }\partial B_R.$$

Set $a=\frac{2n-k}{n-k}$, from $C^0$ estimate, we have
$$(-u^{\varepsilon,R})^{-a}\leq (t^{-1}R)^{\frac{4n-2k}k},\quad\text{on }\partial B_R.$$
So we have 
$$|Du^{\varepsilon,R}|^2(-u^{\varepsilon,R})^{-a}\leq C,\quad\text{on }\partial \Sigma_R,$$ where $C$ is a constant independent of $\varepsilon$ and $R$.

Since
\begin{align}\label{eq718::4}(-u^{\varepsilon,R})^{2-a}\leq \Big((t^{-1}|z|)^{2-\frac {2n}k}\Big)^{-\frac k{n-k}}=t^{-2}|z|^2,\end{align}
and
\begin{align}\label{eq718::5}D\log f^\varepsilon=-(n+1)\frac{\bar z}{|z|^2+\varepsilon^2}.\end{align}
We have
$$(-u^{\varepsilon,R})^{2-a}|D\log f^\varepsilon|^2\leq C_{n}\frac1{t^2}\frac{|z|^4}{(|z|^2+\varepsilon^2)^2}\leq C(n,t).$$
By Theorem \ref{Gradientest}, 
$$ |Du^{\varepsilon,R}|^2(-u^{\varepsilon,R})^{-a}\leq C,$$ where $C$ is independent of $\varepsilon$ and $R$.
Use the $C^0$ estimate once more, we drive that
$$  |Du^{\varepsilon,R}|^2\leq C(-u^{\varepsilon,R})^a
\leq C|z|^{1-\frac{2n}k}.$$

\subsection{Second order estimates}
We will prove the second order estimate of the approximating equations.

\subsubsection{\textbf{The global second order estimate can be reduced to the boundary second order estimate} }
We use the idea of Hou-Ma-Wu \cite{HMW} (see also Chou-Wang \cite{cw2001cpam}) to prove the following estimate.
\begin{theorem}\label{reducesecondtobdry}
Let $u$ be the $k$-subharmonic solution to \eqref{approxeq2} 
and consider $H=u_{\xi\bar \xi}(-u)^{-\frac{n}{n-k}}\psi(P)$. If $(-u)^{-\frac k{n-k}}|D\log f^\varepsilon|^2$ and $(-u)^{-\frac k{n-k}}|D^2\log f^\varepsilon|$ are uniformly bounded which is independent of   $\varepsilon$ and $R$,
then we have
\begin{align}\label{Secondorderestimate-1}
\max_{\Sigma_R}H\leq C+\max_{\partial \Sigma_R}H
\end{align}
where $P=|Du|^2(-u)^{-\frac{2n-k}{n-k}}$, $\psi(t)=(M-t)^{-\sigma}$, $\sigma\leq \frac{a-1}{8a^2}$ and  $M= 2\max\limits_{\Sigma_R}P+1$, $a=\frac{2n-k}{n-k}$, $C$ is a positive constant depending only on $n$, $k$, $\sup\limits_{\Sigma_R} P$, $\sup\limits_{\Sigma_R}(-u)^{-\frac k{n-k}}|D\log f^\varepsilon|^2$ and $\sup\limits_{\Sigma_R}(-u)^{-\frac k{n-k}}|D^2\log f^\varepsilon|$.
\end{theorem}
\begin{theorem}\label{reducesecondtobdry1}
Let $u$ be the $k$-subharmonic solution to \eqref{approxeq2}. Let $\hat w:=-\Big(\frac{\mu^2+|z|^2}{1+\varepsilon^2}\Big)^{1-\frac{n}k}$. Then for sufficient small $\varepsilon$ and $b$, for any unit vector $\xi\in\mathbb R^{2n}$, there holds
$$\max_{\Sigma_R}(\hat w-u-bu_{\xi\xi})\leq \max_{\partial\Sigma_R}(\hat w-u-bu_{\xi\xi}).$$
\end{theorem}

\begin{proof}[\bf Proof of Theorem \ref{reducesecondtobdry}]
For simplicity, we write $f$ instead of $f^\varepsilon$during the proof.\\

Suppose the maximum of $H$ is attained at an interior point $z_0\in \Sigma_R$ along the direction $\xi_0=\frac{\p}{\p z_1}$ .  We can choose the holomorphic coordinate such that  $ \{u_{i\bar j}\}$ is diagonal at $z_0$ and $\lambda_i:=u_{i\bar i}$ with $\lambda_1\ge \lambda_2\ge \cdots\ge\lambda_n$.
The following calculations are  at $z_0$.
 Then we have
$$ 0=\varphi_i=\frac{u_{1\bar 1i}}{u_{1\bar 1}}-(a-1)\frac{u_i}{u}+\sigma\frac{P_i}{M-P}.$$
Denote by $ F^{i\bar j}:=\frac{\partial }{\partial u_{i\bar j}}\log S_k(\partial\bar\partial u)=\frac{S_k^{i\bar j}}{S_k}$, and $ F^{i\bar j,r\bar s}=\frac{\partial ^2 }{\partial u_{i\bar j}u_{r\bar s}}\log S_k(\partial\bar\partial u)=\frac{S^{i\bar j,r\bar s}}{S_k}-\frac{S_k^{i\bar j}S_k^{r\bar s}}{S_k^2}$, $S_k^{i\bar j}:=\frac{\partial}{\partial u_{i\bar j}}S_k(\partial\bar\partial u)$, $S_k^{i\bar j,r\bar s}=\frac{\partial^2}{\partial u_{i\bar j}\partial u_{r\bar s}}S_k(\partial\bar \partial u)$. Then by directly calculation, we have 
\begin{equation}\label{maineq2}\begin{aligned}
0\geq& F^{i\bar j}\varphi_{i\bar j}\\
=&\lambda_1^{-1}F^{i\bar j}u_{1\bar 1i\bar j}-\frac{F^{i\bar i}|u_{1\bar 1i}|^2}{u_{1\bar 1}^2}+(a-1)\frac{F^{i\bar j}u_{i\bar j}}{(-u)}+(a-1)\frac{F^{i\bar i}|u_i|^2}{u^2}+\sigma\frac{F^{i\bar j}P_{i\bar j}}{M-P}+\sigma\frac{F^{i\bar i}|P_i|^2}{(M-P)^2}\\
=&\lambda_1^{-1}F^{i\bar j}u_{1\bar 1i\bar j}-\frac{F^{i\bar i}|u_{1\bar 1i}|^2}{u_{1\bar 1}^2}+\frac{(a-1)k}{(-u)}+(a-1)\frac{F^{i\bar i}|u_i|^2}{u^2}+\sigma\frac{F^{i\bar j}P_{i\bar j}}{M-P}+\sigma\frac{F^{i\bar i}|P_i|^2}{(M-P)^2}\\
:=&\mathrm I+\mathrm{II}+\cdots+\mathrm{VI}.
\end{aligned}\end{equation}
Take the first and second order derivatives to $P$, we have
$$ P_i=|Du|^2_i(-u)^{-a}+|Du|^2((-u)^{-a})_i,$$
and
\begin{equation*}\begin{aligned}
P_{i\bar j}=&|Du|^2_{i\bar j}(-u)^{-a}+|Du|^2_i((-u)^{-a})_{\bar j}+|Du|^2_{\bar j}((-u)^{-a})_i+|Du|^2((-u)^{-a})_{i\bar j}\\
=&\big(u_lu_{\bar li\bar j}+u_{li\bar j}u_{\bar l}+u_{li}u_{\bar l\bar j}+u_{l\bar j}u_{\bar l i}\big)(-u)^{-a}\\
&+a(-u)^{-a-1}\Big(\big(u_lu_{\bar li}+u_{li}u_{\bar l}\big)u_{\bar j}+\big(u_lu_{\bar l\bar j}+u_{l\bar j}u_l\big)u_i\Big)\\
&+a(-u)^{-a-1}|Du|^2u_{i\bar j}+a(a+1)(-u)^{-a-2}|Du|^2u_iu_{\bar j}.
\end{aligned}\end{equation*}
So
\begin{equation}\label{fijpij1}\begin{aligned}
F^{i\bar j}P_{i\bar j}=&F^{i\bar j}\cdot \bigg(\big(u_lu_{\bar li\bar j}+u_{li\bar j}u_{\bar l}+u_{li}u_{\bar l\bar j}+u_{l\bar j}u_{\bar l i}\big)(-u)^{-a}\\
&+a(-u)^{-a-1}\Big(\big(u_lu_{\bar li}+u_{li}u_{\bar l}\big)u_{\bar j}+\big(u_lu_{\bar l\bar j}+u_{l\bar j}u_l\big)u_i\Big)\\
&+a(-u)^{-a-1}|Du|^2u_{i\bar j}+a(a+1)(-u)^{-a-2}|Du|^2u_iu_{\bar j}\bigg)\\
=&2 \mathrm{Re}\{u_l\tilde f_{\bar l}\}(-u)^{-a}+F^{i\bar i}|u_{li}|^2(-u)^{-a}+F^{i\bar i}\lambda_i^2(-u)^{-a}\\
&+2a(-u)^{-a-1}F^{i\bar i}\lambda_i|u_i|^2+2a(-u)^{-a-1}F^{i\bar i}u_{li}u_{\bar l}u_{\bar i}\\
&+ka(-u)^{-a-1}|Du|^2+a(a+1)(-u)^{-a-2}|Du|^2F^{i\bar i}|u_i|^2.
\end{aligned}\end{equation}
and
\begin{equation}\label{fijpij3}\begin{aligned}
F^{i\bar j}P_{i\bar j}
\geq&2 \mathrm{Re}\{u_l\tilde f_{\bar l}\}(-u)^{-a}+ka(-u)^{-a-1}|Du|^2+a(a+1)(-u)^{-a-2}|Du|^2F^{i\bar i}|u_i|^2\\
&+\frac12F^{i\bar i}|u_{li}|^2(-u)^{-a}+\frac12F^{i\bar i}\lambda_i^2(-u)^{-a}\\
&-2a^2(-u)^{-a-2}F^{i\bar i}|u_i|^2|Du|^2-2a^2(-u)^{-a-2}F^{i\bar i}|u_i|^4\\
:=&\mathrm a_1+\cdots+\mathrm a_7.
\end{aligned}\end{equation}
We divide the rest computation into two cases: $\lambda_k\geq \delta\lambda_1$ and $\lambda_k<\delta\lambda_1$.
	
\noindent\textbf{Case 1. } $\lambda_k\geq \delta\lambda_1$. Then
\begin{equation*}\begin{aligned}
\mathrm{II}:=&-\frac{F^{i\bar i}|u_{1\bar 1i}|^2}{u_{1\bar 1}^2}=-F^{i\bar i}\Big|(a-1)\frac{u_i}{u}-\sigma \frac{P_i}{M-P}\Big|^2\\
\geq&-2F^{i\bar i}\bigg((a-1)^2\frac{|u_i|^2}{u^2}+\sigma^2\frac{|P_i|^2}{(M-P)^2}\bigg).
\end{aligned}\end{equation*}
So
\begin{equation*}\begin{aligned}
\mathrm{II}+\mathrm{IV}+\mathrm{VI}:=&-\frac{F^{i\bar i}|u_{1\bar 1i}|^2}{u_{1\bar 1}^2}+(a-1)\frac{F^{i\bar i }|u_i|^2}{u^2}+\sigma \frac{F^{i\bar i}|P_i|^2}{(M-P)^2}\\
\geq&\big((a-1)-2(a-1)^2\big)\frac{F^{i\bar i }|u_i|^2}{u^2}+\big(\sigma-2\sigma^2\big)\frac{F^{i\bar i}|P_i|^2}{(M-P)^2}\\
\geq&\big((a-1)-2(a-1)^2\big)\frac{F^{i\bar i }|u_i|^2}{u^2}.
\end{aligned}\end{equation*}
where the last inequality holds since
$\sigma\leq \frac12$.

By the concavity of $S_k^\frac1k$, we have
$$ \mathrm I:=\lambda_1^{-1}F^{i\bar j}u_{1\bar 1i\bar j}=\lambda_1^{-1}\big((\log f)_{1\bar 1}-F^{i\bar j,r\bar s}u_{i\bar j1}u_{r\bar s\bar 1}\big)\geq \lambda_1^{-1}(\log f)_{1\bar 1}.$$
By \eqref{eq719::1}, we have
\begin{equation}\label{ob1}
F^{i\bar i}\lambda_i^2\geq F^{k\bar k}\lambda_k^2\geq \theta\mathcal F\lambda_k^2\geq \delta^2\theta\mathcal F\lambda_1^2,
\end{equation}
where $\mathcal F=\sum_{i=1}^nF^{i\bar i}$, $\theta=\theta(n,k)$ and we use the assumption of \textbf{Case 1} in the last inequality. Based on \eqref{ob1}, we have the following calculation,
\begin{equation*}\begin{aligned}
\frac14\mathrm a_5+\mathrm a_6+\mathrm a_7:=&\frac18F^{i\bar i}\lambda_i^2(-u)^{-a}-2a^2(-u)^{-a-2}F^{i\bar i}|u_i|^2|Du|^2-2a^2(-u)^{-a-2}F^{i\bar i}|u_i|^4\\
\geq&\frac18F^{i\bar i}\lambda_i^2(-u)^{-a}-4a^2(-u)^{-a-2}\mathcal F|Du|^4\\
\geq& (-u)^{a-2}\mathcal F\bigg(\frac{\delta^2\theta}{8}\big(\lambda_1(-u)^{-a+1}\big)^2-4a^2P^2\bigg)\\
\geq& 0,
\end{aligned}\end{equation*}
where the last inequality holds if we suppose
\begin{equation}\label{Q1} 
\big(\lambda_1(-u)^{-a+1}\big)^2\geq \frac{32a^2}{\delta^2\theta}P^2.
\end{equation}
By Newton-MacLaurin inequality, we have
$$ \frac{S_k}{S_{k-1}}\leq \frac{n-k+1}{nk}S_1.$$
So
$$ \mathcal F=\frac{\sum_{i=1}^n S_{k-1,i}}{S_k}=(n-k+1)\frac{S_{k-1}}{S_k}\geq nkS_{1}^{-1}\geq \frac k{\lambda_1}.$$
Combined with \eqref{ob1}, we have
\begin{equation}\label{ob2}
F^{i\bar i}\lambda_i^2\geq k\delta^2\theta\lambda_1.
\end{equation}
By \eqref{ob2},
\begin{equation*}\begin{aligned}
\frac14\mathrm a_5+\mathrm a_1:=&\frac18F^{i\bar i}\lambda_i^2(-u)^{-a}+2\mathrm{Re}\{u_l\tilde f_l\}(-u)^{-a}\\
\geq& \frac{k\delta^2\theta}{8}\lambda_1(-u)^{-a}-2|Du||D\tilde f|(-u)^{-a}\\
\geq& (-u)^{-1}\bigg(\frac{k\delta^2\theta}{8}\lambda_1(-u)^{-a+1}-2P^{\frac12}|D\tilde f|(-u)^{-\frac{a}{2}+1}\bigg)\\
\geq& 0,
\end{aligned}\end{equation*}
where last inequality holds if we assume 
\begin{equation}\label{Q2}
\lambda_1(-u)^{-a+1}\geq \frac{16}{k\delta^2\theta}|D\tilde f|(-u)^{-\frac a2+1}P^\frac12.
\end{equation}
Note that $ a-1=\frac{n}{n-k}>1$, it follows from \eqref{ob1} that
\begin{equation*}\begin{aligned}
&\frac{\sigma}{M-P}\cdot\frac14\mathrm a_5+\mathrm{II}+\mathrm{IV}+\mathrm{VI}\\
\geq &\frac{\sigma}{M-P}\cdot\frac18F^{i\bar i}\lambda_i^2(-u)^{-a}+\big((a-1)-2(a-1)^2\big)\frac{F^{i\bar i }|u_i|^2}{u^2}\\
\geq &\frac{\sigma}{M-P}\cdot\frac{\delta^2 \theta}8\mathcal F\lambda_1^2(-u)^{-a}+\big((a-1)-2(a-1)^2\big)\mathcal F\frac{|Du|^2}{u^2}\\
=&(-u)^{a-2}\mathcal F\bigg(\frac{\sigma}{M-P}\cdot\frac{\delta^2\theta}{8}\big(\lambda_1(-u)^{-a+1}\big)^2-\big(2(a-1)^2-(a-1)\big)P\bigg)\\
\geq& 0,
\end{aligned}\end{equation*}
where the last inequality holds if we suppose
\begin{equation}\label{Q3}
\big(\lambda_1(-u)^{-a+1}\big)^2\geq \frac{16M}{\sigma\delta^2\theta}\big(2(a-1)^2-(a-1)\big)P.
\end{equation}
By \eqref{ob2},  we have
\begin{equation*}\begin{aligned}
&\frac{\sigma}{M-P}\cdot\frac14\mathrm a_5+\mathrm I\\
\geq& \frac{\sigma}{M-P}\cdot\frac{\delta^2 \theta}8\mathcal F\lambda_1^2(-u)^{-a}-\lambda_1^{-1}(\log f)_{1\bar 1}\\
\geq& \frac{\sigma}{M-P}\cdot\frac{k\delta^2 \theta}8\lambda_1^{-1}(-u)^{-a}-\lambda_1|D^2\log f|\\
=&(-u)^{a-2}\lambda_1^{-1}\bigg(\frac{\sigma}{M-P}\cdot\frac{k\delta^2 \theta}8\big(\lambda_1(-u)^{-a+1}\big)^2-|D^2\log f|(-u)^{-a+2}\bigg)\\
\geq&0,
\end{aligned}\end{equation*}
where the last inequality holds if we suppose 
\begin{equation}\label{Q4}
\big(\lambda_1(-u)^{-a+1}\big)^2\geq \frac{16M}{k\sigma\delta^2\theta}|D^2\log f|(-u)^{-a+2}.
\end{equation}
From assumptions \eqref{Q1}, \eqref{Q2}, \eqref{Q3} and \eqref{Q4}, we have
$$ 0\geq F^{i\bar j}\varphi_{i\bar j}>0.$$ which leads a contradiction. Since $P$, $|D\log f|(-u)^{-\frac a2+1}$ and $|D^2\log f|(-u)^{-a+2}$ are  uniformly bounded, we finish the proof of  {Case 1. }

	\vspace{0.3cm}
\noindent\textbf{Case 2.} $\lambda_k\leq \delta\lambda_1$. By the first order derivatives condition, we have
\begin{equation*}\begin{aligned}
\sigma \sum_{i\geq 2}\frac{F^{i\bar i}|P_i|^2}{(M-P)^2}=&\frac1\sigma \sum_{i\geq 2}F^{i\bar i}\Big|\frac{u_{1\bar 1i}}{u_{1\bar 1}}-(a-1)\frac{u_i}{u}\Big|^2\\
\geq&\frac{\epsilon}\sigma \sum_{i\geq 2}	F^{i\bar i}\Big|\frac{u_{1\bar 1i}}{u_{1\bar 1}}\Big|^2-\frac1\sigma\cdot\frac{\epsilon}{1-\epsilon}(a-1)^2\sum_{i\geq 2}F^{i\bar i}\frac{|u_i|^2}{u^2}.
\end{aligned}\end{equation*}
Putting  the above inequality into \eqref{maineq2}, we have
\begin{equation}\label{maineq3}\begin{aligned}
0\geq& F^{i\bar j}\varphi_{i\bar j}\\
=&\frac{(a-1)k}{(-u)}+\lambda_1^{-1}F^{i\bar j}u_{1\bar 1i\bar j}+\sigma\frac{F^{i\bar j}P_{i\bar j}}{M-P}\\&-\frac{F^{1\bar 1}|u_{1\bar 11}|^2}{u_{1\bar 1}^2}+(a-1)\frac{F^{1\bar 1}|u_1|^2}{u^2}+\sigma \frac{F^{1\bar 1}|P_1|^2}{(M-P)^2}\\
&-\sum_{i\geq 2}\frac{F^{i\bar i}|u_{1\bar 1i}|^2}{u_{1\bar 1}^2}+(a-1)\sum_{i\geq 2}\frac{F^{i\bar i}|u_i|^2}{u^2}+\sigma\sum_{i\geq 2}\frac{F^{i\bar i}|P_i|^2}{(M-P)^2}\\
\ge&\frac{(a-1)k}{(-u)}+\lambda_1^{-1}F^{i\bar j}u_{1\bar 1i\bar j}+\sigma\frac{F^{i\bar j}P_{i\bar j}}{M-P}\\
&-\frac{F^{1\bar 1}|u_{1\bar 11}|^2}{u_{1\bar 1}^2}+(a-1)\frac{F^{1\bar 1}|u_1|^2}{u^2}+\sigma \frac{F^{1\bar 1}|P_1|^2}{(M-P)^2}\\
&-\sum_{i\geq 2}\big(1-\frac\epsilon\sigma\big)\frac{F^{i\bar i}|u_{1\bar 1i}|^2}{u_{1\bar 1}^2}+\Big((a-1)-\frac1{\sigma}\cdot\frac{\epsilon}{1-\epsilon}(a-1)^2\Big)\sum_{i\geq 2}\frac{F^{i\bar i}|u_i|^2}{u^2}\\
:=&\mathrm I'+\mathrm {II}'+\cdots+\mathrm {VIII}'.
\end{aligned}\end{equation}
We take	
\begin{equation}\label{delta1}
\epsilon\leq \min\Big\{\frac14,\frac{3}{8(a-1)}\sigma\Big\},
\end{equation}
then
$$ \mathrm{VIII}'\geq \frac{a-1}{2}\sum_{i\geq 2}\frac{F^{i\bar i}|u_i|^2}{u^2}.$$
Note that 
\begin{equation*}\begin{aligned}
\mathrm{IV}'=&-\frac{F^{1\bar 1}|u_{1\bar 11}|^2}{u_{1\bar 1}^2}=-F^{1\bar 1}\Big|(a-1)\frac{u_1}{u}-\sigma\frac{P_1}{M-P}\Big|^2\\
\geq&-2(a-1)^2F^{1\bar 1}\frac{|u_1|^2}{u^2}-2\sigma^2\frac{|P_1|^2}{(M-P)^2}.
\end{aligned}\end{equation*}
By the choice of $\sigma$, we have
$$ \mathrm{IV}'+\mathrm{V}'+\mathrm{VI}'\geq \big(a-1-2(a-1)^2\big)F^{1\bar 1}\frac{|u_1|^2}{u^2}.$$
Putting the above inequalities into \eqref{maineq3}, 
\begin{equation}\label{maineq4}\begin{aligned}
0\geq &\lambda_1^{-1}F^{i\bar j}u_{1\bar 1i\bar j}-\sum_{i\geq 2}\big(1-\frac\epsilon\sigma\big)\frac{F^{i\bar i}|u_{1\bar 1i}|^2}{u_{1\bar 1}^2}+\frac{(a-1)k}{(-u)}+\sigma\frac{F^{i\bar j}P_{i\bar j}}{M-P}\\
&+\big(\frac{a-1}2-2(a-1)^2\big)F^{1\bar 1}\frac{|u_1|^2}{u^2}+\frac{a-1}{2}\sum_{i\ge 2}\frac{F^{i\bar i}|u_i|^2}{u^2}\\
:=&\mathrm I''+\cdots+\mathrm{VI}''.
\end{aligned}\end{equation}
We have
\begin{equation*}\begin{aligned}
&\mathrm {VI}''+\frac{\sigma}{M-P}(\mathrm a_6+\mathrm a_7)\\:=&\frac{a-1}{2}\sum_{i\ge 2}\frac{F^{i\bar i}|u_i|^2}{u^2}-\frac{\sigma}{M-P}\Big(2a^2(-u)^{-a-2}F^{i\bar i}|u_i|^2|Du|^2+2a^2(-u)^{-a-2}F^{i\bar i}|u_i|^4\Big)\\
=&\sum_{i\ge 2}F^{i\bar i}\frac{|u_i|^2}{u^2}\bigg(\frac{a-1}{2}-\frac{\sigma}{M-P}\cdot 2a^2|Du|^2(-u)^{-a}-\frac{\sigma}{M-P}2a^2|u_i|^2(-u)^{-a}\bigg)\\
&-\frac{\sigma}{M-P}\Big(2a^2(-u)^{-a-2}F^{1\bar 1}|u_1|^2|Du|^2+2a^2(-u)^{-a-2}F^{1\bar 1}|u_1|^4\Big)\\
\geq&\sum_{i\ge 2}F^{i\bar i}\frac{|u_i|^2}{u^2}\bigg(\frac{a-1}{2}-\frac{\sigma P}{M-P}\cdot 4a^2\bigg)-4a^2\frac{\sigma P}{M-P}F^{1\bar 1}\frac{|u_{1}|^2}{u^2}\\
\geq& -(a-1)F^{1\bar 1}\frac{|u_{1}|^2}{u^2},
\end{aligned}\end{equation*}
where the last inequality holds if we take $\sigma\le \frac{a-1}{8a^2}$.

\begin{align}
&\frac{\sigma}{M-P}\big(\frac14\mathrm a_5+\mathrm a_6+\mathrm a_7\big)+\mathrm V''+\mathrm{VI}''\notag\\
\ge &\frac{\sigma}{M-P}\cdot\frac18 F^{i\bar i}\lambda_i^2(-u)^{-a}-\big(\frac{a-1}2+2(a-1)^2\big)F^{1\bar 1}\frac{|u_1|^2}{u^2}\notag\\
\geq& F^{1\bar 1}(-u)^{a-2}\bigg(\frac\sigma{M-P}\cdot\frac18\big(\lambda_1(-u)^{-a+1}\big)^2-(2(a-1)^2+\frac{a-1}2)P\bigg)\notag\\
\geq& 0,\label{08081}
\end{align}
where the last inequality holds if we assume
\begin{equation}\label{Q5}
\big(\lambda_1(-u)^{-a+1}\big)^2\geq \frac{16M}{\sigma }\big(2(a-1)^2+\frac{a-1}2\big)P.
\end{equation}
By
$$ S_k^{-1}(\lambda)S_{k-1}(\lambda|i)\lambda_i^2=S_1(\lambda)-(k+1)\frac{S_{k+1}(\lambda)}{S_k(\lambda)}\geq \frac knS_1(\lambda)\geq \frac kn\lambda_1,$$
we have
\begin{align}
\frac14 \mathrm a_5+\mathrm a_1:=&\frac18F^{i\bar i}\lambda_i^2(-u)^{-a}-2\mathrm{Re}\{u_l\log f_l\}\notag\\
\geq&\frac{k}{8n}\lambda_1(-u)^{-a}-2|Du||D\log f|(-u)^{-a}\notag\\
=&(-u)^{-1}\bigg(\frac{k}{8n}\lambda_1(-u)^{-a+1}-2P^\frac12|D\log f|(-u)^{-\frac a2+1}\bigg)\notag\\
\geq&0,\label{08082}
\end{align}
where the last inequality holds if we assume
\begin{equation}\label{Q6}
\lambda_1(-u)^{-a+1}\geq \frac{16n}{k}P^\frac12|D\log f|(-u)^{-\frac a2+1}.
\end{equation}
By Proposition \ref{lemmaChouWang}, when $ \delta $ is small enough (depending on $\epsilon$ and $\sigma$),
\begin{align}
\mathrm I''+\mathrm {II}'':=&\lambda_1^{-1}F^{i\bar j}u_{1\bar 1i\bar j}-\sum_{i\geq 2}\big(1-\frac\epsilon\sigma\big)\frac{F^{i\bar i}|u_{1\bar 1i}|^2}{u_{1\bar 1}^2}\notag\\
\geq& f^{-1}\lambda_1^{-2}\sum_{i\geq 2}|u_{1\bar 1i}|^2\Big(\lambda_1S_{k-2,1i}-(1-\frac{\epsilon}{\sigma})S_{k-1,i}\Big)+\lambda_1^{-1}(\log f)_{1\bar 1}\notag\\
\geq& -\lambda_1^{-1}|D^2\log f|,\label{08083}
\end{align}
where we use the concavity of $\log S_k$ in the first inequality.

Substituting \eqref{08081}, \eqref{08082} and \eqref{08083} into \eqref{maineq4}, we obtain
\begin{equation}\label{Q7}
0\geq \frac{(a-1)k}{-u}-\lambda_1^{-1}|D^2\log f|.
\end{equation}
Then
\begin{equation}
\lambda_1(-u)^{-a+1}\leq (a-1)k|D^2\log f|(-u)^{-a+2}.
\end{equation}
Since $P$, $|D\log f|(-u)^{-\frac a2+1}$  and $|D^2\log f|(-u)^{-a+2}$ are uniformly bounded, we finish the proof of {Case 2}.	

\end{proof}

\begin{proof}[\bf Proof of Theorem \ref{reducesecondtobdry1}]
Observe that the equation is equivalent to
$$F[u]:=S_k^\frac1k(\partial\bar\partial u)=(f^\varepsilon)^\frac1k.$$ 
Denote by $F^{i\bar j}=\frac{\partial F[u]}{\partial u_{i\bar j}}$ and $F^{i\bar j, k\bar l}=\frac{\partial^2 F[u]}{\partial u_{i\bar j}\partial u_{k\bar l}}$. Now we consider any unit vector $\xi\in\mathbb R^{2n}.$ Differentiating the equation above twice with respect to $\xi$, we obtain
\begin{equation*}\begin{aligned}
F^{i\bar j}u_{\xi\xi i\bar j}=&D_{\xi\xi}F[u]-F^{i\bar j,k\bar l}u_{i\bar j\xi}u_{k\bar l\xi}\geq ((f^\varepsilon)^\frac1k)_{\xi\xi}\\
\geq &-\frac{2(n+k)}k\frac{(f^\varepsilon)^\frac1k}{\varepsilon^2+|z|^2}.
\end{aligned}\end{equation*}
Consider the function 
$$\hat w:=-\big(\frac{\mu^2+|z|^2}{1+\varepsilon^2}\big)^{1-\frac nk}.$$
By the concavity of $S_k^\frac1k$, we have
$$F^{i\bar j}(\hat w_{i\bar j}-u_{i\bar j})\geq F[\hat w]-F[u]=\Big(C_n^k(\frac nk-1)^k(1+\varepsilon^2)^{n-k}\Big)^\frac1k\bigg(\big(\mu^2(|z|^2+\mu^2)^{-k-n}\big)^\frac1k-\big(\varepsilon^2(|z|^2+\varepsilon^2)^{-k-n}\big)^\frac1k\bigg).$$
If $mu>\varepsilon$, then $\frac{\mu^2(|z|^2+\mu^2)^{-k-n}}{\varepsilon^2(|z|^2+\varepsilon^2)^{-k-n}}$ is inscreasing in $|z|$, we have
$$\frac{\mu^2}{\varepsilon^2}\geq\frac{\mu^2(|z|^2+\mu^2)^{-k-n}}{\varepsilon^2(|z|^2+\varepsilon^2)^{-k-n}}\geq \frac{\mu^2(r_0^2+\mu^2)^{-k-n}}{\varepsilon^2(t^2+\varepsilon^2)^{-k-n}}\geq 2^k,\quad\text{in }\Sigma_R,$$
provided $\mu=c_0t$ and $\varepsilon_0\leq 2^{-\frac12}c_0(1+c_0)^\frac{k+n}2$.
So we can take $b\leq \frac{kt^2}{2(n+k)}$ such that 
$$F^{i\bar j}(\hat w-u-au_{\xi\xi})_{i\bar j}\geq 0,\quad\text{in }\Sigma_R.$$
Maximum principle leads that
$$\max_{\overline{\Sigma_R}}(\hat w-u-au_{\xi\xi})\leq \max_{\partial {\Sigma_R}}(\hat w-u-au_{\xi\xi}).$$
\end{proof}
\subsubsection{\textbf{Second order estimate on the boundary $\partial \Sigma_R$}}.

\noindent\textbf{Step1: tangential derivative estimates.}

Consider a point $p\in \partial \Omega$. Without loss of generality, let $ p $ be the origin. Choose the coordinate $z_1,\cdots,z_n$ such that the $x_n$ axis is the inner normal direction to $\partial\Omega$ at $0$. 
Suppose 
$$ t_1=y_1,\ t_2=y_2,\ \cdots, \ t_n=y_n,\ t_{n+1}=x_1,\ t_{n+2}=x_2, \ \cdots ,\ t_{2n}=x_n.$$
Denote by $t'=(t_1,\cdots,t_{2n-1})$. Then around the origin, $\partial\Omega$ can be represented as a graph
$$ t_{2n}=x_n=\varphi(t')=B_{\alpha\beta}t_\alpha t_\beta+O(|t'|^3).$$
Since
$$ u(t',\varphi(t'))=0\quad\text{on }\partial\Omega,$$
we have 
$$ u_{t_\alpha t_\beta}(0)=-u_{t_{2n}}(0)B_{\alpha\beta},\quad \alpha,\beta=1,\cdots,2n-1.$$
It follows that for any $\alpha,\beta=1,\cdots,2n-1$,
\begin{equation}\label{puretangentialfordirichlet} 
|u_{t_\alpha t_\beta}(0)|\leq C,\quad\text{on }\partial\Omega.
\end{equation}

Note that $u\geq \underline u^\varepsilon$ near $\partial\Omega$, $u=\underline u^\varepsilon$ and $0<\underline u^\varepsilon_\nu\leq u_\nu$ on $\partial\Omega$, there exists a smooth function $g$ such that $u=g\underline u$ near $\partial\Omega$, and $g\geq 1$ outside of $\Omega$ nearby $\partial\Omega$. So $\forall\,1\leq i,j\leq n-1$,
\begin{align*}
u_{i\bar j}(0)=g_{i\bar j}(0)\underline u^\varepsilon(0)+g_i(0)\underline u^\varepsilon_{\bar j}(0)+g_{\bar j}(0)\underline u^\varepsilon_i(0)+g(0)\underline u^\varepsilon_{i\bar j}(0).
\end{align*}
Note that $\underline u^\varepsilon=c_0\rho_1$ near $\partial\Omega$, where $\rho_1$ is a given strictly plurisubharmonic function in a neighborhood $\Omega$, $c_0=(1-(1+\frac{s^2}{16+s^2})^{1-\frac{n}k})\tau$, $\tau$ is a constant independent of $\varepsilon$ and $R$ as taken in Lemma \ref{subu0720}. We also have
\begin{align}\label{lowerbdforuijonpo}
S_{k-1}(\{u_{i\bar j}(0)\}_{1\leq i,j\leq n-1})=&c_0^{k-1}g^{k-1}(0)S_{k-1}(\{\rho_{1,i\bar j}(0)\}_{1\leq i,j\leq n-1})\nonumber\\
\geq& c_0^{k-1}g_0^{k-1}C_n^{k-1}(C_n^k)^{\frac{1-k}k}\min_{\partial\Omega}S_k^\frac{k-1}k(\partial\bar \partial \rho_1)>0.
\end{align}

Set $R\geq R_2\geq R_1$, $R_2$ is to be determined later. Consider a harmonic function $h_3$, which is a solution to
\begin{equation}\label{eqforh3}\begin{cases}
\Delta h_3=0&\quad\text{in } B_R\backslash\overline{B_2},\\
h_3=-(1+\varepsilon^2)^{\frac nk-1}(R^2+\varepsilon^2)^{1-\frac nk}&\quad\text{on }\partial B_R,\\
h_3=-t^{\frac{2n}k-2}|z|^{2-\frac {2n}k}&\quad\text{on }\partial B_2.
\end{cases}\end{equation}
Set 
$$\bar h(z):=\tilde h_3(z)=R^{\frac{2n}k-2}h_3(Rz).$$
By maximum principle, we know,
$$\underline{\tilde  u}^\varepsilon\leq \tilde u\leq \bar h,$$
where $\tilde u(z)=R^{\frac{2n}k-2}u(Rz)$.
Note that
$$\bar h=-\big(\frac{1+\frac{\varepsilon^2}{R^2}}{1+\varepsilon^2}\big)^{1-\frac nk},\quad\text{ on }\partial B_1,\quad\text{and}\quad \bar h=-\big(\frac2{Rt}\big)^{2-\frac {2n}k},\quad\text{on }\partial B_{\frac2R}.$$
If we choose $R^2\geq (R_2)^2:=\max\{(R_1)^2,4t^{-2}4(1+\varepsilon_0^2),16\}$, then
$$\bar h\big|_{\partial B_1}\geq \bar h\big|_{\partial B_\frac{2}{R}}.$$
Similarly as in gradient estimates, there is a positive constant $C$, independent of $\varepsilon$ and $R$, such that
$$\bar h_\nu\leq \tilde u_\nu\leq \tilde{\underline u}^{\varepsilon}_\nu\leq C\quad\text{on }\partial B_2.$$
In fact, we can prove that 
$$\bar h_\nu>c_0>0,$$ where $c_0$ is also independent of $\varepsilon$ and $R$.
In fact, we can solve \eqref{eqforh3},
\begin{equation*}\begin{aligned}
\bar  h=&-\frac{-\Big(\frac{1+\frac{\varepsilon^2}{R^2}}{1+\varepsilon^2}\Big)^{1-\frac nk}+\big(\frac{2}{Rt}\big)^{2-\frac nk}}{\big(\frac 2R\big)^{2-N}-1}|z|^{2-N}\\
&-\big(\frac{2}{Rt}\big)^{2-\frac nk}+\Big(-\Big(\frac{1+\frac{\varepsilon^2}{R^2}}{1+\varepsilon^2}\Big)^{1-\frac nk}+\big(\frac{2}{Rt}\big)^{2-\frac nk}\Big)\frac{\big(\frac 2R\big)^{2-N}}{\big(\frac 2R\big)^{2-N}-1}.
\end{aligned}\end{equation*}
Then
\begin{equation*}\begin{aligned}
\bar h_\nu=&\Big(-\Big(\frac{1+\frac{\varepsilon^2}{R^2}}{1+\varepsilon^2}\Big)^{1-\frac nk}+\big(\frac{2}{Rt}\big)^{2-\frac nk}\Big)\big((\frac 2R)^{2-N}-1\big)^{-1}(\frac N2-1)\\
\geq&(\frac N2-1)(2^{\frac nk-1}-1)\big(\frac{1}{1+\varepsilon_0^2}\big)^{1-\frac nk}(2^{N-2}-1)^{-1}\\
>&0.
\end{aligned}\end{equation*}
It follows that there exists a $(\varepsilon,R)$-independent constant $C$, such that
$$C^{-1}R^{1-\frac{2n}k}\leq u_\nu\leq CR^{1-\frac{2n}k}\quad\text{on }\partial B_2,$$ where $\nu$ is the unit outer normal to $\partial B_2$.


For any $p\in \partial B_R$, we choose the coordinate such that $p=(0,\cdots,-R)$. Then near $p$, $\partial B_R$ is locally represented by $t_{2n}=x_n=\varphi(t')=-\sqrt{R^2-\sum_{i=1}^{2n-1}t_i^2}$. Since
$$u(t',\varphi(t'))=-\Big(\frac{R^2+\varepsilon^2}{1+\varepsilon^2}\Big)^{1-\frac nk}\quad\text{on }\partial B_R,$$
we have
$$u_{t_\alpha t_\beta}(p)=-u_{t_{2n}}(p)\frac{\partial^2t_{2n}}{\partial t_\alpha \partial t_\beta}= -R^{-1}u_{t_{2n}}(p)\delta_{\alpha\beta}=R^{-1}u_{\nu}(p)\delta_{\alpha\beta}.	$$
Hence
\begin{align}
&|u_{t_\alpha t_\beta}|\leq CR^{-\frac{2n}k},\quad \alpha,\beta=1,\cdots,2n-1,\label{puretangentialfordirichlet1}\\
&u_{i\bar j}=\frac14(u_{t_{n+i}t_{n+j}}+u_{t_it_j}-iu_{t_it_{n+j}}+iu_{t_{n+i}t_j})\geq CR^{-\frac{2n}k}\delta_{ij},\quad i,j=1,\cdots, n.\label{lowbdforuij}
\end{align}

\noindent\emph{\textbf{Step2: { tangential-normal derivative estimates $\partial \Sigma_R$}}}

Follow the approach by Guan in \cite{Guanbo2014}, we estimate the tangential-normal derivatives on boundary. We first prove the tangential-normal derivatives estimate on $\partial\Omega$. Suppose $0\in \partial\Omega$, to estimate $ u_{t_\alpha t_n}(0) $ for $\alpha=1,\cdots,2n-1$, we consider the auxiliary function 
$$ v=u-\underline u+td-\frac N2 d^2 $$
on $\Omega_\delta=\Omega\cap B_\delta(0)$ with constant $N,t,\delta$ to be determined later. Define a linear operator
$$Lv=F^{i\bar j}v_{i\bar j},$$
where $F^{i\bar j}=\frac{\partial}{\partial u_{i\bar j}}S_k^\frac1k(\partial\bar\partial u) $.
Then 
$$\mathcal{F}=\sum_{i=1}^nF^{i\bar i}=S_k^{\frac1k-1} S_{k-1}(\lambda|i) =(n-k+1)S_k^{\frac1k-1}S_{k-1}\geq C_{n,k}>0.$$

By Lemma \ref{lemma4.5-1},
for $N$ sufficiently large and $t,\delta$ sufficiently small, there holds
\begin{equation*}\begin{cases} 
Lv\leq -\frac{\epsilon}4(1+\mathcal F)&\quad\text{in }\Omega_\delta,\\
v\geq 0&\quad\text{on }\partial\Omega,
\end{cases}\end{equation*}
where $\epsilon>0$ is a uniform constant depending only on subsolution $\underline u$ restricted in a small neighborhood of $\partial\Omega$.

In our setting, $\epsilon$ can be taken independent of $\varepsilon$ and $R$, since $\underline u^\varepsilon=c_0\rho_1$ near $\partial\Omega$, where $\rho_1$ is a given strictly plurisubharmonic function in a neighborhood $\Omega$, $c_0=(1-(1+\frac{s^2}{16+s^2})^{1-\frac{n}k})\tau$, $\tau$ is a constant independent of $\varepsilon$ and $R$ as taken in Lemma \ref{subu0720}.

We use the similar notation as in subsection \ref{subsection3.2}. 
Let
$$\Psi=A_1v+A_2|z|^2-A_3\big((u_{y_n}-\underline u_{y_n})^2+\sum_{l=1}^{n-1}|u_l-\underline u_l|^2).$$
After a similar computation as the boundary tangential-normal derivatives estimate on the pseudoconvex boundary in \ref{subsection3.2}, we see that 
\begin{align*}
L(\Psi\pm T_\alpha (u-\underline u))\leq 0\quad\text{in }\Omega_\delta
\end{align*}
and 
\begin{align*}
\Psi\pm T_\alpha(u-\underline u)\geq 0\quad\text{on }\partial\Omega_\delta,
\end{align*}
when $A_1\gg A_2\gg A_3\gg1$. Therefore
\begin{align}\label{tnestonpo} |u_{t_\alpha x_n}|\leq C\quad\text{on }\partial\Omega.\end{align}

Nextly we prove the tangential-normal derivatives estimate on $\partial B_R$. 
Let 
$$\tilde u(z)=R^{\frac{2n}k-2}u(Rz)\quad\text{and}\quad \tilde {\underline u}^\varepsilon(z)=R^{\frac{2n}k-2}\underline u^\varepsilon(Rz).$$ 
Consider the boundary tangential-normal derivatives estimate on $\partial B_1$. Let $p=(0,\cdots,-1)\in\partial B_1$. Write a defining function $\varrho$ of $B_1$ near $p$ by 
$$\varrho(z)=-x_n-\big(R^2-\sum_{i=1}^{n-1}|z_i|^2-y_n^2\big)^\frac12.$$
Then 
$$|T_\alpha (\tilde u-\tilde {\underline u}^\varepsilon)|\leq C.\quad\text{in  }B_1(0)\cap B_\frac12(p).$$

Let $w=|z|^2-1$, then
$$L(-w)=-\sum_{i=1}^nF^{i\bar i}\leq -C_{n,k}(1+\mathcal F).$$
Let 
$$\Phi=-B_1w+B_2|z-p|^2-B_3\Big(\sum_{l=1}^{n-1}|\tilde {u}_l-\tilde {\underline u}^\varepsilon_l|^2+(\tilde u_{y_n}-\tilde {\underline u}^\varepsilon_{y_n})^2\Big).$$
Similarly, we get 
\begin{align*}
L(\Phi\pm T_\alpha (\tilde u-\tilde {\underline u}^\varepsilon))\leq 0\quad\text{in }\Omega_\delta
\end{align*}
and 
\begin{align*}
\Phi\pm T_\alpha(\tilde u-\tilde{\underline u}^\varepsilon)\geq 0\quad\text{on }\partial\Omega_\delta,
\end{align*}
when $B_1\gg B_2\gg B_3\gg1$.
So  we have
$$|\tilde u_{t_\alpha x_n}|\leq C\quad\text{on }\partial B_1.$$
Therefore
\begin{align}\label{tnestonpb1}| u_{t_\alpha x_n}|\leq CR^{-\frac{2n}k}\quad\text{on }\partial B_R.\end{align}

\noindent\emph{\textbf{Step3: double normal derivative estimates $\partial \Sigma_R$}}\\
By pure tangential derivatives estimate \eqref{puretangentialfordirichlet} and \eqref{puretangentialfordirichlet1}, we have 
$$|u_{y_ny_n}|\leq C\quad\text{on }\partial\Omega \quad\text{and}\quad|u_{y_ny_n}|\leq CR^{-\frac{2n}k}\quad\text{on }\partial B_R.$$
To estimate the double normal deritive $u_{x_nx_n}$, it suffices to estimate $u_{n\bar n}$.
By rotation of $(z_1,\cdots, z_{n-1})$, we may assume that $\{u_{i\bar j}\}_{1\leq i,j\leq n-1}$ is diagonal.
Then
$$ f^\varepsilon=S_k(\partial\bar \partial u)=u_{n\bar n}S_{k-1}(\{u_{i\bar j}\}_{1\leq i,j\leq n-1})+S_k(\{u_{i\bar j}\}_{1\leq i,j\leq n-1})-\sum_{\beta=1}^{n-1}|u_{\beta n}|^2S_{k-2}(\{u_{i\bar j}\}_{1\leq i,j\leq n-1}),$$
It suffices to give a uniform lower positive bound for $S_{k-1}(\{u_{i\bar j}\}_{1\leq i,j\leq n-1})$. 

By \eqref{puretangentialfordirichlet}, \eqref{lowerbdforuijonpo} and \eqref{tnestonpo}, we obtain 
$$u_{n\bar n}(0)\leq C\quad\text{on }\partial\Omega.$$
On the other hand,
$$u_{n\bar n}(0)\geq-\sum_{i=1}^{n-1}u_{i\bar i}\geq -C.$$

By \eqref{puretangentialfordirichlet1}, \eqref{lowbdforuij} and \eqref{tnestonpb1}, we obtain 
\begin{align*}
Cu_{n\bar n}R^{-\frac{2n(k-1)}k}\leq& u_{n\bar n}(0)S_{k-1}(\{u_{i\bar j}\}_{1\leq i,j\leq n-1})\\
=&S_k(\partial\bar \partial u)-S_k(\{u_{i\bar j}(0)\}_{1\leq i,j\leq n-1})+\sum_{\beta=1}^{n-1}|u_{\beta n}(0)|^2S_{k-2}(\{u_{i\bar j}(0)\}_{1\leq i,j\leq n-1})\\
\leq &CR^{-2n}
\end{align*}
Therefore
$$|u_{n\bar n}(0)|\leq CR^{-\frac{2n}k}\quad\text{on }\partial B_R.$$

\noindent\emph{\textbf{Step4: second order derivative estimates in $\Sigma_R$}}\\
As in Theorem \ref{reducesecondtobdry}, let $H=Q(M-P)^\sigma$, $Q=u_{i\bar i}(-u)^{-a+1}$, $P=|Du|^2(-u)^{-a}$. Suppose the maximum of $H$ is obtain at a boundary point $z_0\in\partial\Sigma_R$. Then
\begin{align}\label{eq718::3}Q=(M-P)^\sigma H\leq M^\sigma H(z_0)\leq M^\sigma Q(z_0)(M-\max_{\Sigma_R} P)^{-\sigma}=\max_{\partial \Sigma_R}Q(M-\max_{\Sigma_R}P)^{-\sigma}.\end{align}
Note that $P$ is bounded (uniformly in $\varepsilon$ and $R$). By \eqref{eq718::4} and \eqref{eq718::5}, $$(-u)^{2-a}|D\log f^\varepsilon|^2\leq C(n,k,t)\quad\text{and}\quad
(-u)^{2-a}|D^2\log f^\varepsilon|^2\leq C(n,k,t).$$
By Theorem \ref{reducesecondtobdry}, if the maximum of $H$ is obtained at a interior point, there is a positive constant $C$ independent of $\varepsilon$ and $R$ such that 
$Q\leq C.$
Combined with \eqref{eq718::3}, there is a positive constant $C$ independent of $\varepsilon$ and $R$ such that 
$$Q\leq C \quad\text{in }\Sigma_R.$$
Then we obtain
\begin{align}\Delta u\leq C (-u)^{a-1}\leq C|z|^{-\frac {2n}k} \quad\text{in }\Sigma_R.\end{align}

By boundary second order derivative estimates and $C^0$ estimate, we obtain that for any unit vector $\xi\in\mathbb R^{2n}$,
$$\max_{\partial {\Sigma_R}}(\hat w-u-au_{\xi\xi})\leq C.$$Hence
$$ u_{\xi\xi}\leq C,\quad\text{in }\overline\Sigma_R.$$
$u$ is subharmonic since $u$ is $k$-admissible,  then
$$-C\leq u_{\xi\xi}\leq C \quad\text{in }\Sigma_R.$$ In conclude, we get
\begin{align} |D^2u|\leq C,\quad\text{in }\Sigma_R.\end{align}

\section{Proof of Theorem \ref{main07201}}\label{section5}
\subsection{Uniqueness}
The uniqueness follows from the comparison principle for $k$-subharmonic solutions of the complex $k$-Hessian equation in bounded domains in Lemma \ref{comparison0718} by Blocki \cite{BlockiAIF}. 

Suppose $u$ and $v$ are two solutions to \eqref{case1Equa1.1}. For any $z_0\in\mathbb C^n\setminus\Omega$, there exists $R_0$ such that $z_0\in B_{R_0}(0)\backslash\Omega$. Since $u(z)\rightarrow 0$, $v(z)\rightarrow 0$ as $|z|\rightarrow \infty$, $\forall\,\varepsilon>0$, there exists $R\gg R_0$ such that 
$$v-\varepsilon\leq u\leq v+\varepsilon\quad\text{in }\mathbb C^n\setminus B_R.$$ By the comparison principle Lemma \ref{comparison0718}, 
$$v-\varepsilon\leq u\leq v+\varepsilon\quad\text{in } B_R\setminus \Omega.$$ Note that $z_0\in B_{R_0}(0)\setminus\Omega\subset B_R(0)\setminus\Omega$, we have
$$v(z_0)-\varepsilon\leq u(z_0)\leq v(z_0)+\varepsilon.$$ Let $\varepsilon\rightarrow 0$, we obtain that
$u(z_0)=v(z_0)$. Since $z_0$ is arbitrary, $u=v$ in $\mathbb C^n\backslash \Omega$.

\subsection{The existence and $C^{1,1}$-estimates}
The existence follows from the uniform $C^2$-estimates for $u^{\varepsilon, R}$. The proof is similar as that in \cite{gb2007imrn} by Guan.

For any fixed $M_0>R_2$, for the solution to \eqref{approxeq2}, by the $C^2$ estimates, we have
$$\|u^{\varepsilon,R}\|_{C^2(\overline \Sigma_{M_0})}\leq C_1\quad\text{ independent of }\varepsilon, R \text{ and }M_0.$$
for all $R\geq M_0$. By the Evans-Krylov theory, we obtain for $0<\alpha<1$,
$$\|u^{\varepsilon,R}\|_{C^{2,\alpha}(\overline \Sigma_{M_0})}\leq C_2(\varepsilon, M_0)\quad\text{ independent of }R .$$
By compactness, we can find a sequence $R_j\rightarrow\infty$ such that
$$u^{\varepsilon,R_j}\rightarrow u^\varepsilon\quad\text{in }C^2(\overline \Sigma_{M_0}),$$
where $u^\varepsilon$ satisfies
\begin{equation*}\label{approxeq1}
\begin{cases}
H_k(u^\varepsilon)=f^\varepsilon&\quad\text{in }\Sigma_{M_0},\\
u=-1&\quad\text{on }\partial\Omega,
\end{cases}
\end{equation*}
and 
\begin{align*}
	&C^{-1}|z|^{2-\frac{2n}k}\leq -u^{\varepsilon}(z)\leq C|z|^{2-\frac{2n}k},\\
	&|Du^{\varepsilon}(z)|\leq C|z|^{1-\frac{2n}k},\\
	&|\p\bar\p  u^{\varepsilon}(z)|\leq C|z|^{-\frac{2n}k},\\
	&|D^2u^{\varepsilon}(z)|\leq C.
\end{align*}
Moreover,
\begin{equation*}
	\|u^{\varepsilon}\|_{C^{2,\alpha}(\overline \Sigma_{M_0})}\leq C_2(\varepsilon, M_0)\quad \text{for any }\ M_0>R_2.
	\end{equation*}
By the classical Schauder theory,  $u^{\varepsilon}$ is smooth.

 By the above decay estimates for $u^{\varepsilon}$, for any sequence $\varepsilon_j\rightarrow 0$, there is a subsequence of $\{u^{\varepsilon_j}\}$ converging to a function $u$ in $C^{1,\alpha}$ norm on any compact subset of $\mathbb C^n\setminus \Omega$ ( for any $0<\alpha<1$). Thus $u\in C^{1,1}(\mathbb C^n\setminus \Omega)$ and satisfies the disired estimates \eqref{decay10720}. By the convergence theorem of the complex $k$-Hessian operator proved by Trudinger-Zhang in \cite{TZ2014} (see also Lu \cite{Lu2015}), $u$ is a  solution to \eqref{case1Equa1.1}.

{\bf Acknowledgements:}
The second author was supported by  National Natural Science Foundation of China (grants 11721101 and 12141105) and National Key Research and Development Project (grants SQ2020YFA070080). The third author was supported by NSFC grant  No. 11901102.

\bibliographystyle{plain}
\bibliography{2022-08-08-Gao-Ma-Zhang}

\end{document}